\documentclass[final,5p,times,twocolumn]{elsarticle}

\usepackage{graphicx}
\usepackage{amssymb}
\usepackage{amsthm}
\usepackage{lipsum}
\usepackage{xcolor}
\usepackage{comment}
\usepackage{mathtools}
\usepackage{gensymb}
\usepackage{float}
\usepackage{dblfloatfix}
\usepackage{graphicx}
\usepackage{multicol} 
\usepackage{commath}
\usepackage{framed}
\usepackage{nomencl}
\usepackage{siunitx}
\makenomenclature

\renewcommand*\nompreamble{\begin{multicols}{2}}
\renewcommand*\nompostamble{\end{multicols}}
\newcommand{\nomunit}[1]{%
\renewcommand{\nomentryend}{\hspace*{\fill}#1}}

\usepackage{caption}
\usepackage{subcaption}
\usepackage[ruled,vlined]{algorithm2e}
\SetKwProg{Init}{Initialization:}{}{}
\theoremstyle{definition}

\DeclareCaptionLabelSeparator{emdash}{\textemdash}

\newcommand{\tr}{{\scalebox{.63}{$\mathsf{T}$}}} 

\newcommand{\cmm}[1]{{\color{magenta}#1}}

\iftrue 
\usepackage[colorlinks=true]{hyperref}
	\usepackage{xcolor}
	\hypersetup{
		colorlinks,
		linkcolor={blue!55!black},
		citecolor={red!55!black},
		urlcolor={blue!80!black}
	}
	\allowdisplaybreaks
	
\fi

\journal{Energy and Buildings}
\begin{document}
\begin{frontmatter}
\title{Distributed Model Predictive Control of Buildings and Energy Hubs}

\address[label1]{Department of Mechanical Engineering, Swiss Federal Institute of Technology, Lausanne}
\address[label2]{Automatic Control Laboratory, Swiss Federal Institute of Technology, Z\"{u}rich}
\address[label3]{Empa, Urban Energy Systems Laboratory, Überlandstrasse 129, 8600 Dübendorf, Switzerland}

\author[label1,label2]{Nicolas~Lefebure\, }
\author[label2]{Mohammad~Khosravi\, \corref{cor1}}
\author[label2]{Mathias~Hudoba~de~Badyn\,  \corref{cor1}}
\author[label2,label3]{Felix~B\"{u}nning\, \corref{cor1}}
\author[label2]{John~Lygeros\, }
\author[label1]{Colin~Jones\, }
\author[label2]{Roy~S.~Smith\, }
\address{nicolas.lefebure@epfl.ch,  \{khosravm, mbadyn\}@control.ee.ethz.ch, felix.buenning@empa.ch,
colin.jones@epfl.ch,
\{lygeros, rsmith\}@control.ee.ethz.ch
}
\cortext[cor1]{These authors contributed equally to the work.}
\begin{abstract}
Model predictive control (MPC) strategies can be applied to the coordination of energy hubs to reduce their energy consumption.
Despite the effectiveness of these techniques, their potential for energy savings are potentially underutilized due to the fact that energy demands are often assumed to be fixed quantities rather than controlled dynamic variables. 
The joint optimization of energy hubs and buildings' energy management systems can result in higher energy savings. 
This paper investigates how different MPC strategies perform on energy management systems in buildings and energy hubs. 
We first discuss two MPC approaches; centralized and decentralized. 
While the centralized control strategy offers optimal performance, its implementation is computationally prohibitive and raises privacy concerns. 
On the other hand, the decentralized control approach, which offers ease of implementation, displays significantly lower performance. 
We propose a third strategy, distributed control based on dual decomposition, which has the advantages of both approaches. 
Numerical case studies and comparisons demonstrate that the performance of distributed control is close to the performance of the centralized case, while maintaining a significantly lower computational burden, especially in large-scale scenarios with many agents. 
Finally, we validate and verify the reliability of the proposed method through an experiment on a full-scale energy hub system in the NEST demonstrator in D\"{u}bendorf, Switzerland.
\end{abstract}
\end{frontmatter}
\section{Introduction}

\begin{table*}[ht]   
\begin{framed}
\mbox{}

\nomenclature{\parbox[t]{1.5cm}{$N_{\text{HP}}$, $N_{\text{HB}}$, $N_{\text{S}}$, $N_{\text{NC}}$}}{\parbox[t]{5.75cm}{Number of heat pumps, electric boilers, tanks, uncontrolled buildings. In this case $N_{\text{S}}$ $\geq$ $N_{\text{HP}}$, $N_{\text{HB}}$, $N_{\text{NC}}$.} \nomunit{[-]}}
\nomenclature{$N_{\text{C}}$}{Number of controlled buildings. \nomunit{[-]}} 
\nomenclature{$x^s_i(k)$}{Average temperature of tank $i$ at instant $k$. \nomunit{[$\degree$C]}} 
\nomenclature{$\Delta b^s_i(k)$}{Heat balance in tank $i$ at instant $k$. \nomunit{[$kW$]}} 
\nomenclature{$d^s_i(k)$}{\parbox[t]{5.75cm}{Heat demand by the buildings at instant $k$ connected to tank $i$.} \nomunit{[$kW$]}} 
\nomenclature{$z^s_i(k)$}{\parbox[t]{5.75cm}{Binary number determining if $u^s_i(k)$ is switched on or off.} \nomunit{[-]}} 
\nomenclature{$u^s_i(k)$}{\parbox[t]{5.75cm}{Heat supplied to tank $i$ by heat pump and electric boiler at instant $k$.} \nomunit{[$kW$]}} 
\nomenclature{$x^c_i(k)$}{\parbox[t]{5.75cm}{Vector containing room temperatures in controlled building $i$ at instant $k$.} \nomunit{[$\degree$C]}} 
\nomenclature{$u^c_i(k)$}{\parbox[t]{5.75cm}{Vector containing heat supplied in rooms of controlled building $i$ at instant $k$.} \nomunit{[$kW$]}} 
\nomenclature{$d^c_i(k)$}{\parbox[t]{5.75cm}{Vector of solar irradiance and ambient temperature in controlled building $i$ at instant $k$.} \nomunit{[$kW$]}} 
\nomenclature{$z^c_i(k)$}{\parbox[t]{5.75cm}{Binary number determining if $u^c_i(k)$ is switched on or off.} \nomunit{[-]}} 
\nomenclature{$u^{s_{\text{net}}}_{i}(k)$}{\parbox[t]{5.75cm}{Stacked vector containing the heat demand of connected controlled buildings to tank $i$.} \nomunit{[$kW$]}} 
\nomenclature{$u^{c_{\text{net}}}_{i}(k)$}{\parbox[t]{5.75cm}{Stacked vector containing the heat taken from each connected tank to controlled building $i$.} \nomunit{[$kW$]}} 
\nomenclature{$\mathcal{L}_{i}$}{\parbox[t]{5.75cm}{Set containing the indexes of all controlled buildings connected to tank $i$.} \nomunit{[-]}} 

\printnomenclature[0.6in] 
\end{framed}
\end{table*}

Recent developments have considerably diversified and expanded the technologies used to harvest and manage energy. 
The utilization of all these technologies to operate multiple buildings effectively and cooperatively  has led to the concept of energy hubs~\cite{darivianakis2015stochastic,darivianakis2017power, 4042137, MURRAY20194204}. 
An energy hub comprises different energy production, conversion, and storage capabilities, whose objective is to efficiently manage energy resources to handle time-varying production/consumption mismatches~\cite{Bayod}. 
With the arrival of renewable energy sources and the proliferation of prosumers, energy-producing environments will use an increasing number of energy technologies. 
The concept of the energy hub is therefore promising for the efficient management of such environments.

Despite their interconnected nature, energy hubs have no impact on the effective consumption of energy by their connected consumers, such as buildings. 
According to~\cite{SwissFed:2018}, the Swiss Federal Office of Energy (SFOE) states that in 2021, Swiss buildings will consume approximately 100 TWh. 
This corresponds to 45\% of the total energy demand nationwide. 
On a global scale, buildings consume 32\% of the world's total energy demand~\cite{URGEVORSATZ201585}. 
The SFOE declared that the Energy Strategy 2050 (Switzerland’s new energy policy) aims to reduce the energy consumption of Swiss buildings to 55 TWh by 2050. 
To achieve this objective, buildings need to be considered as active participants with energy management systems -- working hand in hand with energy hubs to accommodate and manage new energy technologies to reduce the world's energy consumption.

Energy hubs can benefit from advanced control methods, such as model predictive control (MPC,~\cite{MORARI1999667}), to provide stable and accurate energy management strategies that are in accordance with the (as of now, uncoordinated) exterior energy supply/demand of buildings. 
Two contrasting control strategies are commonly present: decentralized and centralized. 
The principle of the decentralized approach is that control is local and that there is no communication between the energy hub and the different local controllers of the buildings. 
For example,~\cite{7483144} proposes a decentralized MPC approach for the management of energy hubs. 
In the centralized approach, a single master controller is designed to compute the control actions of all the controllers of the energy hub and the buildings, accounting explicitly for interactions between the energy hub and the buildings.
In~\cite{5275230}, a centralized MPC to control multi-energy systems is demonstrated in simulation.

The same advanced control techniques are also particularly suitable for building control, since by employing  a predictive approach, the evolution of building behavior, weather fluctuations, variation in energy prices, and so forth,~can all be taken into account to optimally adapt the heating/cooling control policy to reduce energy consumption while maintaining thermal comfort constraints. As MPC approaches are often limited by the need for fast computing power, it is particularly adaptable to building automation, since the slow timescales of the thermodynamic processes facilitate implementation in real-time. The idea of using MPC for building automation is not new, as there have already been many attempts both in simulation \cite{7536001, MA201292, 8814615, 5759140} and on real systems~\cite{MA20141282, HILLIARD2017326, 7087366}. 
A comprehensive survey on MPC for building automation is provided by~\cite{drgovna2020all}.

The approaches in the literature have not yet considered buildings as controllable entities in combination with energy hubs, even though the thermal mass of buildings can effectively act as energy storage~\cite{DOMINKOVIC2018949}. 
With data-driven models of building dynamics becoming increasingly available~\cite{SMARRA20181252, BUNNING2020109792, khosravi2019data, MADDALENA2020104211, KATHIRGAMANATHAN2021110120, bunning}, building models can be included in the energy hub control problem with moderate effort. Furthermore, in this configuration, the building controllers provide the computation, communications and control framework.
Therefore, instead of deploying MPC strategies for energy hubs that treat buildings simply as energy demands, one can envision simultaneously controlling both energy hubs and buildings, leading to more significant energy savings. 

Therefore, this work aims to consider energy hubs and buildings as cooperatively controlled entities with local constraints. These considerations raise questions and concerns regarding the privacy of the building unit occupants and the computational feasibility of such an approach.
We address this issue by using a distributed control system, where the controllers of the units perform calculations separately from each other and communicate by only sharing virtual prices.
The introduced scheme is compared to a centralized and a decentralized predictive control approach through extensive numerical experiments modeling a physical system. 
Finally, the proposed method is tested on an experimental configuration equipped with an energy hub and a building in
D\"ubendorf, Z\"urich. 

The paper is structured as follows. 
The modeling environments of the energy hubs and the buildings are defined in Section~\ref{problem_statement}. 
The centralized, decentralized, and distributed control structures are described in Section~\ref{methodology}. 
All three controllers are compared via numerical simulations in Section~\ref{Results}. 
Experimental results are presented in Section~\ref{experiment}, and the paper is concluded in Section~\ref{sec:conclusion}.

\section{Problem Statement}
\label{problem_statement}

\label{models}
\begin{figure}[t!]
    \centering
    \includegraphics[width=\columnwidth]{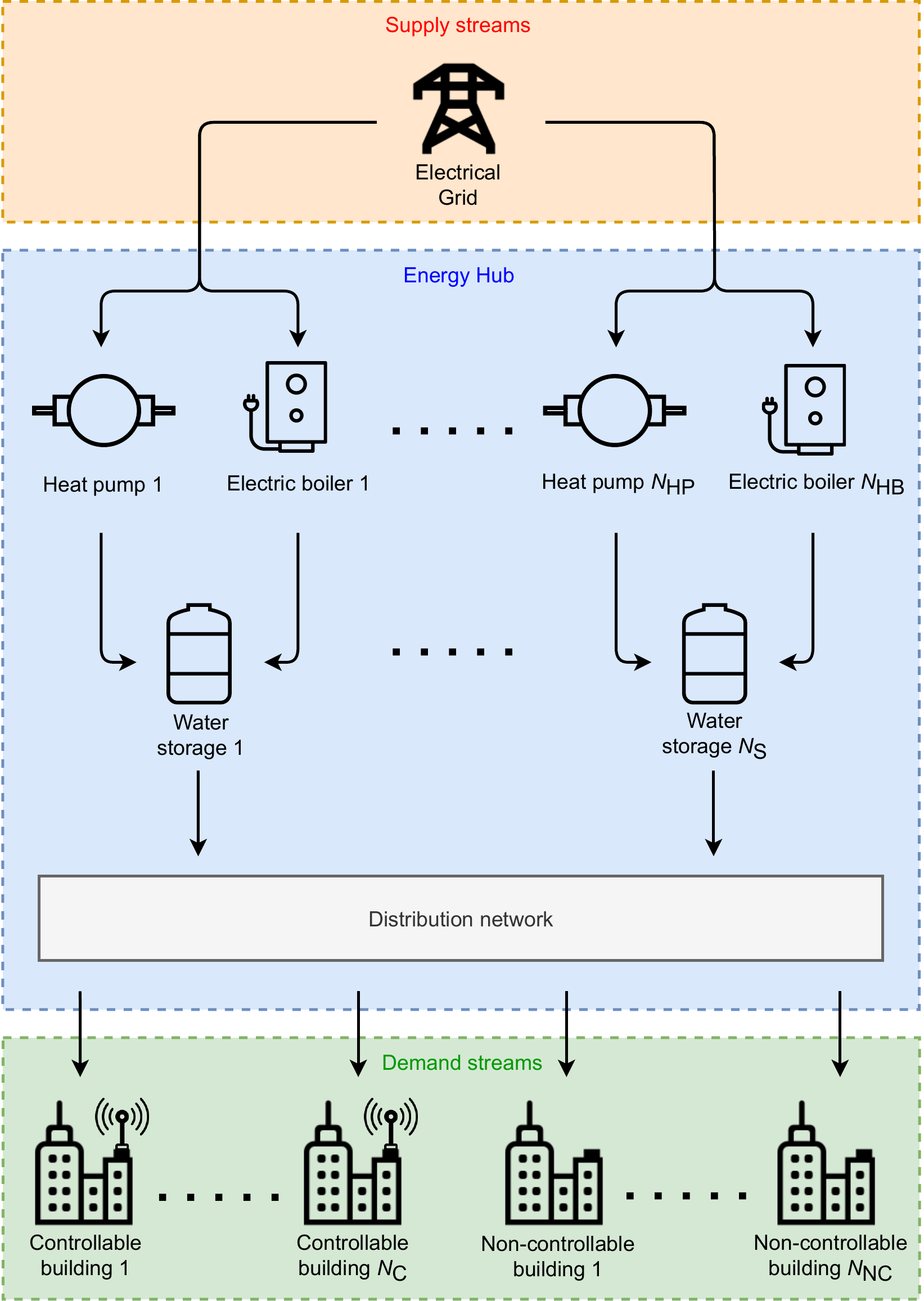}
    \caption{The configuration of the energy hub and buildings.}
    \label{arch1}
\end{figure}

We first define the architecture of the environment in which the energy hubs and the buildings interact. The environment is made of infrastructures and technologies that were readily available for performing experiments. To demonstrate the potential benefits of distributed control, we focus only on serving a heating demand of buildings served by systems made of tanks, heat pumps, and boilers that receive electricity from the grid. The methods could extend to more general settings, as the technologies that constitute an energy hub can be very diverse~\cite{SADEGHI2019114071}, but would obscure our main point, which is the relative benefits between decentralized, centralized, and distributed control. Accordingly, in this study, energy hubs consist of storage units, like water tanks for thermal energy storage, conversion units such as heat pumps and electric boilers, and the network units like heat distribution networks. These are connected to external energy supply and demand, e.g., the electrical grid and the heating demands of buildings. Figure~\ref{arch1} illustrates the architecture scheme of the energy hub and the exterior streams. 
The hub, supply streams, and demand streams are shown by dashed rectangles respectively in blue, red, and green color.
The connection arrows depict the energy flows between the several structures. Starting from the top, the supply stream consists of the local electrical grid supplier. For ease of discussion, it is assumed that the supply stream is constant. Accordingly, we unify all the different potential energy sources, such as hydropower and solar power, as one single source, shown here by the electrical grid supplier.  In the middle, the energy hub comprises two levels. On the first level, conversion components  consist of $N_{\text{HP}}$ heat pumps and $N_{\text{HB}}$ electric boilers, which are connected to the electrical grid.  The next level consists of $N_{\text{S}}$ water storage tanks. Every tank is supplied by a maximum of one heat pump and one electric boiler, and so $N_{\text{HP}} \geq N_{\text{HB}}, N_{\text{S}}$. At the bottom of Figure~\ref{arch1}, one can see that the tanks are serving the heating demand streams via a heat distribution network. This network defines the existing links between each building and each storage tank. We assume that $N_{\text{C}}$ of the buildings are controlled, i.e., they coordinate the control inputs with the energy hub. Meanwhile, the rest of the buildings, $N_{\text{NC}}$, are uncontrolled, i.e., they do not coordinate with the energy hub, and their demand is rather seen as a disturbance. Lastly, we assume that all uncontrolled buildings connected to one tank can be lumped into a single demand, so the number of uncontrolled buildings is less than to the number of water tanks, i.e., $N_{\text{S}} \geq N_{\text{NC}}$.

Next, we introduce the models for each of the components. 
For water tanks, we consider a discrete-time linear equation to describe the dynamics. 
The equation that captures the dynamical behavior of the $i$th water tank in the energy hub at time $k$ is given by
\begin{equation}
    x^s_i(k+1) = A^s_ix^s_i(k) + B^s_i\Delta b^s_i(k) + E^s_id^s_i(k),
    \label{stor}
\end{equation}
where $x^s_i(k)$ denotes the state variable which is the average temperature of the tank, $d^s_i(k)$ is the disturbance and defined as the heating demand by the uncontrolled buildings, and $\Delta b^s_i(k)$ denotes the input balance which is specified with more details later in this section. Note that the superscript ``$s$'' on each variable stands for ``supplier''.
The state and input variables are limited by polytopic operational constraints as follows
\begin{align}
    H^s_ix^s_i(k)  &\leq h^s_i + \epsilon^s_i(k),
    \label{ineq1}\\
    G^s_iu^s_i(k) &\leq g^s_i(1 - z^s_i(k)) + \tilde{g}^s_i z^s_i(k),
    \label{ineq2}
\end{align}
where $h^s_i$ is the vector denoting state constraints, $g^s_i$ and $\tilde{g^s_i}$ are vectors for input constraints, and $H^s_i,G^s_i$ are matrices of appropriate dimensions. 
Note that $z^s_i(k)\in\{0,1\}$ is a binary variable that determines whether the input is switched on or off. 
We emphasize that the variable $z^s_i(k)$ is necessary only when the input $u^s_i(k)$ is constrained by a non-zero lower bound.
Hence, we enforce that $g^s_i$ and $\tilde{g^s_i}$ are complementary vectors, i.e., if an input command needs to be switched off, $z^s_i(k)$ can switched between 0 and 1 and $\tilde{g^s_i}$ is non-zero, otherwise $z^s_i(k)$ is defined as null and $g^s_i$ is non-zero. Furthermore, similar to~\cite{6730917}, the slack variable $\epsilon^s_i(k)$ is introduced to relax the state constraints and guarantee feasibility at each time-step.  
Since the states with bound constraints are the temperatures, it is preferred to penalize possible constraint violation rather than the algorithm terminates due to infeasibility.

The heat supplied to the water tanks is provided by heat pumps and electric boilers. 
It is assumed that the coefficient of performance of the heat pumps is higher than that of the boilers. 
Their operating range has a non-zero lower capacity limit ($u_{\text{HP}} \in [u_{\text{low}}, u_{\text{high}}]$, with $u_{\text{low}}, u_{\text{high}} > 0$). The electric boiler has a larger range of operation and thus acts as a backup if the heat pump is not able to satisfy the energy demand under these constraints. 

The dynamics of the controlled buildings are also described by a discrete-time linear time-invariant system. 
More precisely, for the $j$th building connected to the hub, at instant $k$, we have
\begin{align}
    x^c_j(k+1) = A^c_jx^c_j(k) + B^c_ju^c_j(k) + E^c_jd^c_j(k),
    \label{apart}
\end{align}
where $x^c_j(k)$ is the vector of state variables which are the temperatures of the zones in the building,  $u^c_j(k)$ denotes the vector of heating powers in the zones, and $d^c_j(k)$ is the vector of disturbances including solar radiation and ambient temperature. Note that the superscript ``$c$'' over each variable stands for ``consumer''.
For the states and inputs, we have the following polytopic operational constraints 
\begin{align}
    H^c_jx^c_j(k)  & \leq h^c_j  + \epsilon^c_j(k),
    \label{ineq3} \\
    G^c_ju^c_j(k) & \leq g^c_j(1 - z^c_j(k)) + \tilde{g}^c_j z^c_j(k),
    \label{ineq4}
\end{align}
where $h^c_j$ is the vector for state constraints, $g^c_j$ and $\tilde{g}^c_j$ are the vectors for input constraints, and 
$H^c_j,G^c_j$ are matrices of appropriate dimensions. 
The variable $z^c_j(k)\in\{0,1\}$ is binary, and determines whether the input is switched on or off. The slack variable $\epsilon^c_j(k)$ is used here to soften the state constraints.

In contrast to controlled buildings, the energy consumption in uncontrolled buildings  is determined by local users or unknown control systems. 
From the perspective of the suppliers, the demand of such buildings is represented as a disturbance in equation~\eqref{stor}. 
The model of the uncontrolled buildings is therefore reduced to a heating demand forecast that can be obtained using historical data and available techniques in machine learning~\cite{khosravi2017automated}. 
We rely on a feed-forward ANN (Artificial Neural Network)~\cite{ABIODUN2018e00938} to provide day-ahead forecasts based on ambient conditions and time features. 
Additionally, the predictions are improved with two correction methods, one based on the forecasting error auto-correlation, and one based on online learning. 
The methodology is described in~\cite{BUNNING2020109821} and validated  in~\cite{Bnning2020MachineLA} where an experiment is performed  employing this forecast technique and robust MPC for the purpose of frequency regulation.

The heat balance between energy suppliers and consumers must be satisfied at each time-step $k$ in accordance with the assumed distribution network. We model this requirement through the following equality constraints
\begin{align}
    &\Delta b^s_i(k)  = C^{s}_i u^{s}_i(k) + B^{s_{\text{net}}}_i u^{s_{\text{net}}}_{i}(k), 
    \quad 
    i = 1, \dots, N_\text{S},
    \label{bal1} \\
    &u^c_j(k)  = B^{c_{\text{net}}}_j u^{c_{\text{net}}}_{j}(k), 
    \qquad\qquad\quad\ \ \  
    i = 1, \dots, N_\text{C},
    \label{bal2} \\
    &C^{s_{\text{net}}}_i u^{s_{\text{net}}}_{i}(k)  = \sum_{j \in \mathcal{L}_i} C^{c_{\text{net}}}_{ij} u^{c_{\text{net}}}_{j}(k), 
    \quad\ \ \  
    i = 1, \dots, N_\text{S}
    \label{bal3}
\end{align}
where equation~\eqref{bal1} describes the input-output energy balance for each storage tank.
More precisely, for storage tank $i$, the term $u^{s}_i(k)$ refers to the stream incoming to the storage which includes the heat power supplied by the heat pump and the electric boiler, and, the second term $u^{s_{\text{net}}}_{i}(k)$ is a vector containing the heat demand of the connected controlled buildings. 
Moreover, equation~\eqref{bal2} defines that the input vector $u^c_j(k)$ of each controlled building must be equal to a vector $u^{c_{\text{net}}}_{j}(k)$ containing the amount of heat taken from each connected tank. 
Finally, equation~\eqref{bal3} ensures the balance between the heat quantity taken from each tank  and the heat consumed by each connected controlled building. 
For tank $i$, the set $\mathcal{L}_{i}$ contains the indices of all controlled buildings connected to that tank.

\section{Control methodologies}
\label{methodology}
In this section, we use the previously defined dynamics and  constraints to derive a suitable control strategy for determining the control action of the energy hub. To this end, we present different forms of MPC control structures for the energy hub environment. More precisely, centralized, decentralized and distributed approaches are respectively presented in Section~\ref{centralizedMPC}, Section~\ref{DecentralizedMPC} and Section~\ref{DistributedMPC}, respectively. 
In Section~\ref{SolvDistributedMPC}, an algorithm for the deployment of the distributed control structure is proposed.

\subsection{Centralized MPC}
\label{centralizedMPC}
Given the models, the constraints, the current measurements, and the disturbance forecasts, the centralised MPC controller computes an optimal control input by solving an open-loop optimization problem formulated over a prediction horizon of length $N$. Then, the first timestep of the computed control sequence is applied to the plant~\cite{8767123}. 
Repeating this procedure and shifting the horizon produces a closed-loop controller as new measurements are integrated into the optimisation problem at each step.  
 
We have described two classes of agents, namely suppliers $s$ and consumers $c$. 
A supplier agent corresponds to a water storage supplied by a heat pump and an electric boiler, which is disturbed by a building heating demand. 
A consumer agent is a controlled building.
The cost function\footnote{We choose a quadratic cost function as it provides a balance between minimizing total costs and peaks.} of the $i$th agent is defined as
\begin{equation}
    f^\alpha_i(v^\alpha_i) =  v_{i}^{\alpha}{}^\tr  S^\alpha_{i} v^\alpha_{i},
\end{equation}
where $\alpha \in \{s,c\}$ depending on whether $i$ is a supplier or a consumer, and where we define the decision variable and cost matrix respectively as
\begin{equation}\label{eqn:vialpha_Sialpha}
v^\alpha_{i}
=
\begin{bmatrix}
u^\alpha_i \\
\epsilon^\alpha_i \\
x^\alpha_i \\
z^\alpha_i \\
u^{\alpha_{\text{net}}}_i
\end{bmatrix}, \qquad S^\alpha_{i}=\begin{bmatrix}
Q^\alpha_{i} & 0 & 0 & 0 & 0\\
0 & R^\alpha_{i} & 0 & 0& 0\\
0 & 0 & 0 & 0 & 0\\
0 & 0 & 0 & 0 & 0 \\
0 & 0 & 0 & 0 & 0 
\end{bmatrix},
\end{equation}
where
 $Q^\alpha_i$ and $R^\alpha_i$ are user-defined weighting matrices, and $u^\alpha_{i}$,  
$\epsilon^\alpha_{i}$, 
$x^\alpha_{i}$, 
$z^\alpha_{i}$
and $u^{\alpha_{\text{net}}}_i$ are  vectors of stacked quantities respectively for input, slack, state, binary and network variables, as follows
\begin{equation*}
\begin{split}
u^\alpha_{i} &=\begin{bmatrix}
u^\alpha_{i}(0) \\
\vdots \\
u^\alpha_{i}(N-1)
\end{bmatrix}, \; \epsilon^\alpha_{i}=\begin{bmatrix}
\epsilon^\alpha_{i}(1) \\
\vdots \\
\epsilon^\alpha_{i}(N)
\end{bmatrix}, \;
x^\alpha_{i}=\begin{bmatrix}
x_{i}^\alpha(1) \\
\vdots \\
x_{i}^\alpha(N)
\end{bmatrix},
\\
z^\alpha_{i}&=\begin{bmatrix}
z_{i}^\alpha(0) \\
\vdots \\
z_{i}^\alpha(N-1)
\end{bmatrix}, \;
u^{\alpha_{\text{net}}}_i=\begin{bmatrix}
u^{\alpha_{\text{net}}}_i(1) \\
\vdots \\
u^{\alpha_{\text{net}}}_i(N)
\end{bmatrix}.
\end{split}
\end{equation*}

In this notation, the inequality constraints in equations~\eqref{ineq1} and~\eqref{ineq2} as well as~\eqref{ineq3} and~\eqref{ineq4} can be written as
\begin{equation}
\begin{bmatrix}
\tilde{G}^\alpha_i  & 0 & 0 & g^\alpha_i - \tilde{g}^\alpha_i & 0\\
0 & -\tilde{I} & \tilde{H}^\alpha_i & 0 & 0 \\
0 & -\tilde{I} & 0 & 0 & 0
\end{bmatrix}\begin{bmatrix}
u^\alpha_{i} \\
\epsilon_{i} \\
x^\alpha_{i} \\
z^\alpha_{i} \\
u^{\alpha_{\text{net}}}_i
\end{bmatrix} \leq \begin{bmatrix}
g^\alpha_{i} \\
h^\alpha_{i} \\
0
\end{bmatrix},\label{eq:ineqverbose}
\end{equation}
where $\tilde{G}^\alpha_i = I_{N} \otimes G^\alpha_i$, $\tilde{g}^\alpha_i = I_{N} \otimes g^\alpha_i$, $\tilde{H}^\alpha_i = I_{N} \otimes H^\alpha_i$, $ \tilde{I} = I_{N} \otimes I$ and $\otimes$ denotes the Kronecker product. 
The inequality \eqref{eq:ineqverbose} can be written in the standard form
\begin{equation}
     \mathcal{G}^\alpha_{i}(v^\alpha_i) \leq 0.
\end{equation}
where $\mathcal{G}^\alpha_{i}(v^\alpha_i)$ is an affine function of $v^\alpha_i$. 

We can write the equality constraints in equations \eqref{stor}, \eqref{apart}, \eqref{bal1}, \eqref{bal2} and~\eqref{bal3} respectively as
\begin{align}
    &
    \left(I - \tilde{A}^s_{i}\right) x^s_{i}-\tilde{B}^s_{i} u^s_{i} - \tilde{B}^{s_{\text{net}}}_i u^{s_{\text{net}}}_i  = \tilde{E}^s_{i} d^s_{i}+c^{s,x_0}_{i}, 
    \\&
    \left(I - \tilde{A}^c_{j}\right) x^c_{j}-\tilde{B}^c_{j} u^c_{j} = \tilde{E}^c_{j} d^c_{j}+c^{c,x_0}_{j}, 
    \\&
    u^c_{j}-\tilde{B}^{c_{\text{net}}}_j u^{c_{\text{net}}}_j  = 0,
    \\&
    \tilde{C}_i^{s_{\text{net}}} u^{s_{\text{net}}}_i  - \sum_{j \in \mathcal{L}_i} \tilde{C}_{ij}^{c_{\text{net}}} u^{c_{\text{\text{net}}}}_j  = 0,
\end{align}
where the matrices $ \tilde{A}^\alpha_{i}$ and vectors $c^{\alpha,x_0}_{i}$ are defined as
\begin{equation*}
    \tilde{A}^\alpha_{i}=\begin{bmatrix}
0 & 0 & \cdots & & 0 \\
A^\alpha_{i} & 0 & \cdots & & 0 \\
0 & A^\alpha_{i} & & & 0 \\
\vdots & & \ddots & & \vdots \\
0 & \cdots & & A^\alpha_{i} & 0
\end{bmatrix}, \;
c^{\alpha,x_0}_{i} =\begin{bmatrix}
A^\alpha_{i} x^\alpha_{i}(0) \\
0 \\
\vdots \\
0
\end{bmatrix},
\end{equation*}
and where $\tilde{B}^\alpha_{i}=I_{N} \otimes B^\alpha_{i}C^\alpha_{i}$, $\tilde{B}^{\alpha_{\text{net}}}_{i}=I_{N} \otimes B^\alpha_{i}B^{s_{\text{net}}}_{i}$, $\tilde{E}^\alpha_{i}=I_{N} \otimes E^\alpha_{i}$ and $\tilde{C}^{\alpha_{\text{net}}}_{ij}=I_{N} \otimes C^{\alpha_{\text{net}}}_{ij}$. 
In a similar way, the equality constraints can be written in standard form
\begin{align}
   &\mathcal{F}^s_{i}(v^s_{i}, (u^{c_{\text{net}}}_{j})_{j \in \mathcal{L}_{i}})  = 0,\\
   &\mathcal{F}^c_{j}(v^c_{j})  = 0,
\end{align}
where $\mathcal{F}^s_{i}(v^s_{i}, (u^{c_{\text{net}}}_{j})_{j \in \mathcal{L}_{i}})$ and $\mathcal{F}^c_{j}(v^c_{j})$ are affine functions of $v^s_{i}$, $v^c_{j}$ and $(u^{c_{\text{net}}}_{j})_{j \in \mathcal{L}_{i}})$.
The resulting centralized optimization problem is as follows:
\begin{equation}\tag{P1}\!\!
    \begin{array}{cll}
\underset{\substack{v^s_1,\dots, v^s_{N_\text{S}} \\ v^c_1,\dots, v^c_{N_\text{C}}}}{\text{min}} & \sum\limits_{i=1}^{N_\text{S}} f^s_i(v^s_i) + \sum\limits_{j=1}^{N_\text{C}} f^c_j(v^c_j) \\
\text { s.t. } &  \mathcal{F}^s_{i}(v^s_{i},  (u^{c_{\text{net}}}_{j})_{j \in \mathcal{L}_{i}}) = 0, & i = 1,\dots,N_\text{S}, \\
&  \mathcal{G}^s_{i}(v^s_i) \leq 0,  & i = 1,\dots,N_\text{S}, \\
& \mathcal{F}^c_{j}(v^c_{j}) = 0,    & j = 1, \dots, N_\text{C},\\
& \mathcal{G}^c_{j}(v^c_{j}) \leq 0, & j = 1, \dots, N_\text{C}. 
\end{array}
\label{P1}
\end{equation}
Here, the supplier’s equality constraint $\mathcal{F}^s_{i}$ (in particular, the heat balance equations in ~\eqref{bal3}) depends on the consumer’s control input $(u^{c_{\text{net}}}_{j})_{j \in \mathcal{L}_{i}}$, which couples their dynamics.
The resulting problem is a Mixed Integer Quadratic Program (MIQP,~\cite{hijazi2017convex}), due to the binary variables in the constraints.
\subsection{Decentralized MPC}
\label{DecentralizedMPC}
The decentralized MPC approach is simply a partitioning of the centralized approach~\cite{christofides2013distributed}. 
The control problem is divided into $N_\text{S} + N_\text{C}$ local problems of smaller size. 
In this context, agents define their own optimization problems and make control decisions independently from each other. 
These decisions rely exclusively on local information (such as measurements, forecasts or control decisions) and there is no negotiation between agents during the optimization process. 
More precisely, for each $i$, the $i$th agent decides on its decision variable $v^\alpha_i$ such that its own cost function $f^\alpha_i(v^\alpha_i)$ is minimized while its proper constraints $\mathcal{F}^\alpha_{i}$ and $\mathcal{G}^\alpha_{i}$ are satisfied. 
From the point of view of the agent $i$, other decision variables $v^\alpha_j$ with $j \neq i$ are ignored, i.e. the energy balance constraints coupling them together are applied in the physical system, but not taken into account in the local decision problem.

\subsection{Distributed MPC}
\label{DistributedMPC}
In the distributed MPC approach~\cite{christofides2013distributed}, agents pass information to one another to facilitate solving a coupled, global optimization problem.
In such a setting it is often desired to limit the communication between agents, due to, for example, to computational complexity, and privacy concerns.

In the MPC setting considered here, one approach to designing distributed controllers is based on the dual-decomposition method~\cite{10.5555/993483}. 
First, the centralized problem is decomposed into agent-based sub-problems. 
The sub-problems are then driven towards the global optimal solution of the centralized problem the \emph{dual problem} that acts as a coordinator between the sub-problems through a shared dual variable. 
Our approach to this problem is inspired by \cite{5160224}. 
We start by introducing variables $\mathrm{r}_j^{(i)}$,  $i=1,\ldots,N_{\text{S}}$ with $j\in\mathcal{L}_i$. These introduced variables contain local versions of the coupled variable $u^{c_{\text{net}}}_{j}$. So problem~\eqref{P1} can be written as
\begin{equation}\tag{P2}\!\!
    \begin{array}{cll}
\underset{\substack{v^s_1,\dots, v^s_{N_\text{S}} \\ v^c_1,\dots, v^c_{N_\text{C}}}}{\text{min}} & \sum\limits_{i=1}^{N_\text{S}} f^s_i(v^s_i) + \sum\limits_{i=1}^{N_\text{C}} f^c_i(v^c_i) \\
\text { s.t. } &  
\mathcal{F}^s_{i}(v^s_{i}, \mathrm{r}^{}_{i} ) = 0, & i = 1,\dots, N_\text{S}, \\
&  \mathcal{G}^s_{i}(v^s_i) \leq 0, & i = 1,\dots,N_\text{S}, \\
& \mathcal{F}^c_{j}(v^c_{j}) = 0, & j = 1, \dots, N_\text{C},\\
& \mathcal{G}^c_{j}(v^c_{j}) \leq 0, & j = 1, \dots, N_\text{C}, \\
& 
(u^{c_{\text{net}}}_{j})_{j \in \mathcal{L}_{i}} - \mathrm{r}_i = 0, & i = 1, \dots, N_\text{S},\\
\end{array}
\label{P2}
\end{equation}
where $\mathrm{r}_i$ is defined as $\mathrm{r}_i=(r^{(i)}_{j})_{j \in \mathcal{L}_{i}}$, for each $i=1,\ldots,N_{\text{S}}$. 
The bottom equality constraint ensures that the introduced local variables $\mathrm{r}_i$ are equal to the shared, coupled variable $(u^{c_{\text{net}}}_{j})_{j \in \mathcal{L}_{i}}$.
The Lagrangian for problem~\eqref{P2} then becomes
\begin{equation}
    \begin{split}
&\mathcal{L}(\bar{v}^s, \bar{v}^c, \bar{\lambda}^s, \bar{\mu}^s, \bar{\lambda}^c, \bar{\mu}^c, \bar{p}, \bar{r}) 
\\& \ \ =    
\sum\limits_{i=1}^{N_\text{S}} f^s_i(v^s_i) + \sum\limits_{j=1}^{N_\text{C}} f^c_j(v^c_i) 
+ \sum\limits_{i=1}^{N_\text{S}} \lambda^{s}_i{}^\tr \mathcal{F}^{s}_i(v^s_{i}, \mathbf{r}_{i} )
\\&  \ \  \ \ 
 + \sum\limits_{i=1}^{N_\text{S}} \mu^{s}_i{}^\tr \mathcal{G}^{s}_i(v^s_i)
 + \sum\limits_{i=1}^{N_\text{C}} \lambda^{c}_j{}^\tr \mathcal{F}^c_{j}(v^c_{j}) 
+ \sum\limits_{j=1}^{N_\text{C}} \mu^{c}_j{}^\tr \mathcal{G}^c_{j}(v^c_{j}) 
\\&  \ \  \ \ 
+ \sum\limits_{i=1}^{N_\text{S}} \mathrm{p}_i^{\tr} ( \mathrm{u}^{c_{\text{net}}}_i- \mathrm{r}_{i}),
\end{split}
\end{equation}
where, for each $i=1,\ldots,N_{\text{S}}$, 
we have
$\mathrm{u}^{c_{\text{net}}}_i= (u^{c_{\text{net}}}_{j})_{j \in \mathcal{L}_{i}}$.
The vectors $\bar{v}^s$, $\bar{\lambda}^s $, $\bar{\mu}^s$, $\bar{r}$, $\bar{v}^c$, $\bar{\lambda}^c$, $\bar{\mu}^c$ and $\bar{p}$ are respectively defined as
\begin{equation*}
\begin{split}
\bar{v}^s & =\begin{bmatrix}
v^s_1 \\
\vdots \\
v^s_{N_\text{S}}
\end{bmatrix}, \;
    \bar{\lambda}^s=\begin{bmatrix}
\lambda^s_1 \\
\vdots \\
\lambda^s_{N_\text{S}}
\end{bmatrix}, \; \bar{\mu}^s=\begin{bmatrix}
\mu^s_1 \\
\vdots \\
\mu^s_{N_\text{S}}
\end{bmatrix}, \;
\bar{r}=\begin{bmatrix}
\mathrm{r}_1 \\
\vdots \\
\mathrm{r}_{N_\text{S}}
\end{bmatrix}, \;
\\
 \bar{v}^c & =\begin{bmatrix}
v^c_1 \\
\vdots \\
v^c_{N_\text{C}}
\end{bmatrix}, \;
\bar{\lambda}^c =\begin{bmatrix}
\lambda^c_1 \\
\vdots \\
\lambda^c_{N_\text{C}}
\end{bmatrix}, \;
\bar{\mu}^c=\begin{bmatrix}
\mu^c_1 \\
\vdots \\
\mu^c_{N_\text{C}}
\end{bmatrix}, \; \bar{p}=\begin{bmatrix}
\mathrm{p}_1 \\
\vdots \\
\mathrm{p}_{N_\text{C}}
\end{bmatrix}. \;
\end{split}
\end{equation*}
For $i = 1, \dots, N_\text{S}$, $\lambda^s_{i}$ and $\mu^s_{i}$ are  the local Lagrange multipliers associated to the supplier $i$; 
for $j = 1, \dots, N_\text{C}$, $\lambda^c_{j}$ and $\mu^c_{j}$ are  the local Lagrange multipliers associated to the consumer $j$; and 
for $i=1,\ldots,N_{\text{S}}$, the vectors $\mathrm{p}_i = (p_j^{(i)})_{j \in \mathcal{L}_{i}}$ are the global dual variables shared between the suppliers and the consumers. 
Note that the problem is now separable and so we can write the following dual decomposition:
\begin{equation}
\begin{split}
&
\max\limits_{\substack{\bar{\lambda}^s, \; \bar{\mu}^s, \; \bar{\lambda}^c, \; \bar{\mu}^c, \; \bar{p}}} 
\  
\min\limits_{\substack{\bar{v}^s, \; \bar{v}^c, \; \bar{r}}} \mathcal{L}(\bar{v}^s, \bar{v}^c, \bar{\lambda}^s, \bar{\mu}^s, \bar{\lambda}^c, \bar{\mu}^c, \bar{p}, \bar{r}) 
\\&  \ \ =
\max\limits_{\bar{p}}
\max\limits_{\substack{\bar{\lambda}^s, \; \bar{\mu}^s, \; \bar{\lambda}^c, \; \bar{\mu}^c}}  
\min\limits_{\substack{\bar{v}^s, \; \bar{v}^c, \; \bar{r}}}
\Bigg\{\sum\limits_{i=1}^{N_\text{S}} f^s_i(v^s_i) 
+ \sum\limits_{i=1}^{N_\text{S}} \lambda^{s}_i{}^\tr \mathcal{F}^{s}_i(v^s_{i},  \mathrm{r}_i)
\\& \ \  \ \ \ \ 
 + \sum\limits_{i=1}^{N_\text{S}} \mu^{s}_i{}^\tr \mathcal{G}^{s}_i(v^s_i) + \sum\limits_{j=1}^{N_\text{C}} f^c_j(v^c_j) 
\\& \ \  \ \  \ \
+ \sum\limits_{j=1}^{N_\text{C}} \lambda^{c}_j{}^\tr \mathcal{F}^c_{j}(v^c_{j}) + \sum\limits_{j=1}^{N_\text{C}} \mu^{c}_j{}^\tr \mathcal{G}^c_{j}(v^c_{j})
\\&  \ \  \ \  \ \
+ \sum\limits_{i=1}^{N_\text{S}} \mathrm{p}_i^{\tr} (\mathrm{u}^{c_{\text{net}}}_{i} - \mathrm{r}_i)\Bigg\}
\\&  \ \ =
\max\limits_{\bar{p}}
\Bigg\{\sum\limits_{i=1}^{N_\text{S}}  \max\limits_{\lambda^s_{i}, \; \mu^s_{i}} \min\limits_{v^s_{i}, \mathrm{r}_i}\bigg\{f^s_{i}(v^s_{i})
 +
\lambda^{s} _{i}{}^\tr \mathcal{F}^s_{i}(v^s_{i}, \mathrm{r}_i) 
 \\&  \ \  \ \  \ \ 
  + \mu^{s}_{i}{}^\tr \mathcal{G}^s_{i}(v^s_{i}) - \sum_{j \in \mathcal{L}_i} p_j^{(i)}{}^{\tr} 
  r_j^{(i)}\bigg\} + \sum\limits_{j=1}^{N_\text{C}} \max\limits_{\lambda^c_{j}, \; \mu^c_{j}} \min\limits_{v^c_{j}}\bigg\{f^c_{j}(v^c_{j})
\\&  \ \  \ \ \ \ 
  + \lambda^{c}_{j}{}^\tr \mathcal{F}^c_{j}(v^c_{j}) + \mu^{c}_{j}{}^\tr \mathcal{G}^c_{j}(v^c_{j}) + p^{(j)}{}^{\tr}  u^{c_{\text{net}}}_{j}\bigg\}\Bigg\},
\end{split}\!\!\!\!\!\!\!\!\!\!
\label{Dual_Decomp}
\end{equation}
where vector $p^{(j)}$ is defined as
\begin{equation*}
    p^{(j)} = \sum\limits_{\{i|j\in\mathcal{L}_i\}}p_i^{(j)}, \qquad j=1,\ldots,N_\text{C},
\end{equation*}
and we have used the fact that
\begin{equation}
    \sum\limits_{i=1}^{N_\text{S}} \mathrm{p}_i^{\tr}  \mathrm{u}_i^{c_{\text{net}}} 
    =
    \sum\limits_{j=1}^{N_\text{C}}  p^{(j)}{}^{\tr}  u^{c_{\text{net}}}_{j}.
\end{equation}
We assume that the Slater's condition hold for problem~\eqref{P2}, which is further discussed later.
Subsequently, the strong duality holds \cite{10.5555/993483}. Accordingly, from the dual decomposition in equation~\eqref{Dual_Decomp}, one can distribute the problem  across the agents with one problem per agent based on the dual ascent method \cite{10.5555/993483}. 

More precisely,  we introduce the following iterative scheme:
\begin{equation*}\tag{P3}
\begin{array}{l}
\begin{aligned}
\begin{pmatrix} v_{i}^{s,+} \\ \mathrm{r}_{i}^{+} \end{pmatrix} := & \; \underset{v^s_i,\ \mathrm{r}_i}{\text{argmin}}
& & f^{s}_i (v^{s}_i) + \sum_{j \in \mathcal{L}_i} \mathcal{R}_{ij}(r_j^{(i)}) - p_j^{(i)}{}^{\tr} 
r_j^{(i)} \\
& \quad\; \text{s.t.}
& & \mathcal{F}^s_{i}(v^s_{i}, r_i) = 0,\\
& \; &&   \mathcal{G}^{s}_i(v^s_i) \leq 0.\\
\end{aligned}\\
\\
\begin{aligned}
\begin{pmatrix}v_{j}^{c,+}\end{pmatrix} := & \; \underset{v^c_{j}}{\text{argmin}}
& &  f^c_{j} (v^c_{j}) + p^{(j)}{}^{\tr}  u^{c_{\text{net}}}_{j}\\
& \quad\; \text{s.t.}
& & \mathcal{F}^c_{j}(v^c_{j}) = 0,\\
& \; &&   \mathcal{G}^c_{j}(v^c_{j}) \leq 0, \\
\end{aligned}  \\
\\
\begin{aligned}
\begin{pmatrix}\mathrm{p}_{i}^{+}\end{pmatrix} & \; :=  & & \mathrm{p}_{i}+\kappa(\mathrm{u}_{i}^{c_{\mathrm{net}}}-\mathrm{r}_{i}).\\
\end{aligned}
\end{array}\!\!\!\!\!\!\!\!\!\!
\label{P3}
\end{equation*}
where 
$\kappa$ is the step size, the superscript ``$+$'' denotes the iteration update, and   
$\mathcal{R}_{ij}(r_j^{(i)}) := \rho_{ij} \|r_j^{(i)}\|^2$
is a \emph{regularization term}. 
Given real positive scalars $\rho_{ij}$, for $i = 1, \dots, N_\text{S}$ and $j\in\mathcal{L}_i$, these regularization terms are employed to improve the convergence of the proposed scheme.
Towards the same goal of regularizing the problem, we may replace $S_i^{\alpha}$, introduced in \eqref{eqn:vialpha_Sialpha} with $S_i^{\alpha}+\delta I$, where $I$ is the identity matrix and $\delta$ is a small positive weight.  

The distribution among the agents means that the original centralized problem~\eqref{P1} splits into $N_\text{S} + N_\text{C}$ separate optimization problems that can be solved iteratively in parallel as described in problem~\eqref{P3}.
Following this, for $i=1,\ldots,N_{\text{S}}$ the price vectors $\mathrm{p}_i$ are updated based on the last equation in problem~\eqref{P3}. 
Using the introduced scheme, the decision variables converge to an optimal point and the dual variables converge to an optimal dual point \cite{10.5555/993483}. 

\subsection{Implementation of the distributed MPC problem}
\label{SolvDistributedMPC}
\begin{figure}[h!]
    \centering
    \includegraphics[width=\columnwidth]{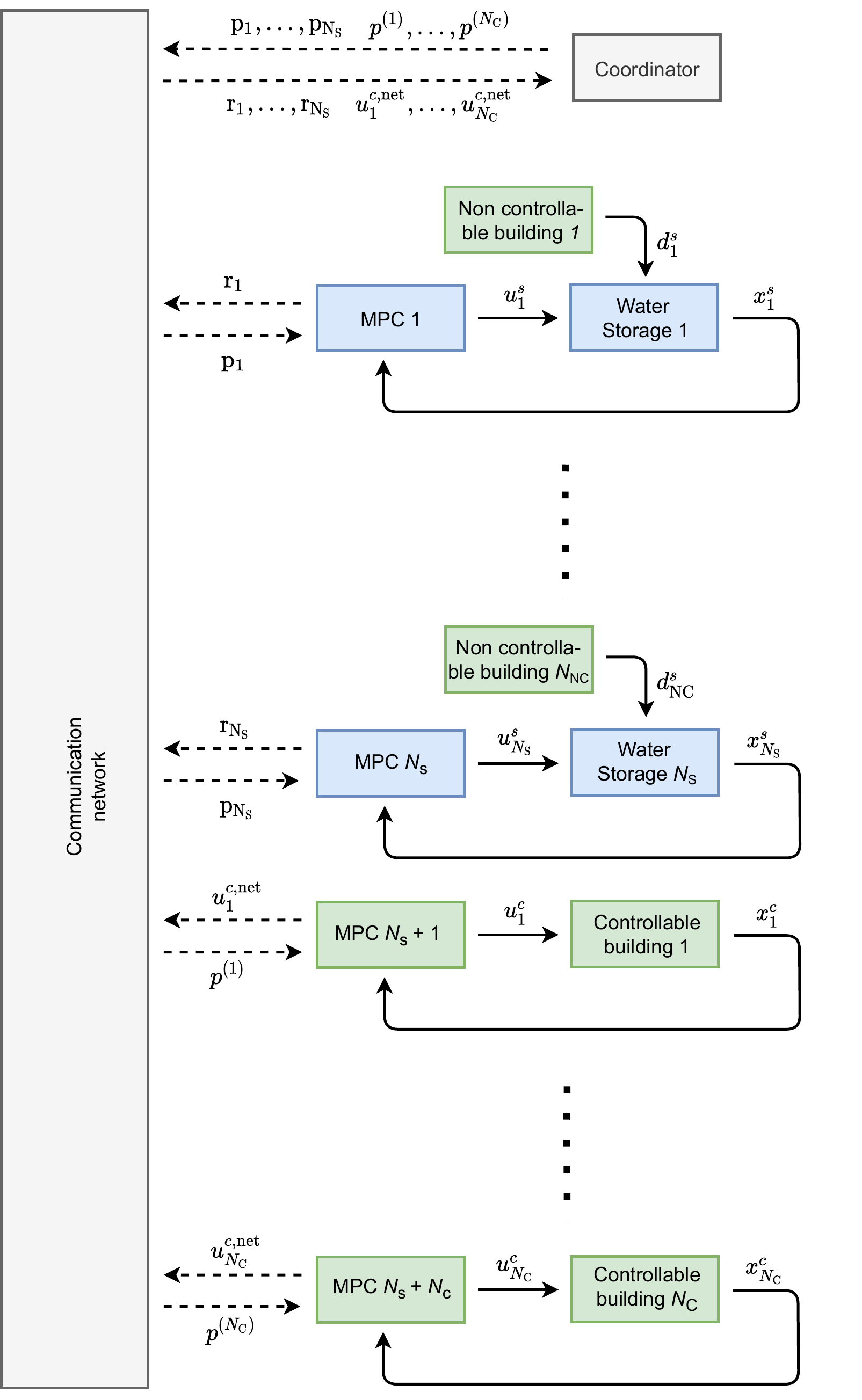}
    \caption{Distributed MPC communication architecture (solid arrows indicate the actions/measurements applied at each time-step while dashed arrows indicate the information shared iteratively between each time-step).}
    \label{DMPC_diagram}
\end{figure}

The procedure to find the solution to the distributed problem~\eqref{P3} can be formulated as follows. 
Given an initial dual variable $p_i^{(j)}$ for each link between a consumer $j$ and a supplier $i$, the minimization problems in~\eqref{P3} are solved by their respective agents. 
Agents $i$ and $j$ then share the optimal values found for $\mathrm{r}_i$ and $u^{c_{\text{net}}}_{j}$ with an external unit through a shared communication network. 
Then, the external unit computes the update of the prices in~\eqref{P3} using the sub-gradient method (see \cite{10.5555/993483} for details), and broadcasts them to the corresponding agents.
Alternately, one of the agents can also compute the price update and broadcast it across the network to the corresponding agent.
The agents recalculate their optimal values and the process is repeated until a convergence criterion is achieved. 
Figure~\ref{DMPC_diagram} represents the architecture of the communication occurring between the agents.

It should be noted that when solving problem \eqref{P2} with the sub-gradient method, the convergence towards the solution can be slow as the local cost functions are quadratic and the dual cost is linear. 
Accordingly, we improve the convergence rate by including a small quadratic regularization term in the cost function of each supplier. 
This term is $\mathcal{R}_{ij}(r_j^{(i)})$, as shown in problem~\eqref{P3}.

Due to Slater's condition and strong duality, solving problem~\eqref{P3} is equivalent to solving the regularized version of problem~\eqref{P1}. 
However, in this paper the variables $z^\alpha_{i}$ belong to a discrete set and thus compromise the convexity of both problems. 
As a result, a duality gap can appear between the dual problem~\eqref{P3} and the primal~\eqref{P1}. 
In order to resolve this duality gap, we relax the problem~\cite{d2003relaxations}, i.e. replace some of the binary variables $z^\alpha_{i}$ (the choice of which is discussed below) with continuous variables in the interval $[0,1]$. 

The solution of the relaxed problem $R$(\ref{P3}) then matches that of the corresponding relaxed centralized problem $R$(\ref{P1}). 
However, the solution will be sub-optimal when projected back to the original problem with binary variables. In order to address this issue, we propose a two-step method~\cite{ding2014two}: 
\begin{enumerate}
    \item Solve the relaxed problem $R$(\ref{P3}) until convergence, i.e., we reach threshold $\epsilon_{\text{tol,r}}$ in the variation of cost function. 
    \item Fix the binary variables $z_i^{\alpha}$ (based on the solution of the relaxed problem, i.e., when it is above a defined threshold value $z_{\text{bound}}$, it is set to $1$, and otherwise, it is set to $0$), and then, solve problem~\eqref{P3} until convergence $\epsilon_{\text{tol}}$ is achieved, i.e., the variation of cost function is less than a given threshold $\epsilon_{\text{tol}}$.
\end{enumerate}

A known MPC practice is move-blocking \citep{CAGIENARD2007563}, where decision variables at the end of the prediction horizon are constrained to be equal, as these have a small effect on the optimality of the implemented control input at the current time step. With the same reasoning, we only relax the variables after a certain horizon $N_{\text{relax}}$~\cite{flamm2019twostage} and keep the first $N_{\text{relax}}$ ones as binaries.
This way we can improve the estimation of the first $N_{\text{relax}}$ binary variables~\cite{mehanna2014feasible}, while lowering the complexity of the problem compared to the non-relaxed MIQP.
If $N_{\text{relax}}\geq 2$, we refer to this problem as the \emph{semi-relaxed} MIQP.

\begin{figure}[t!]

\centering
     \begin{subfigure}[b]{\columnwidth}
         \centering
         \includegraphics[width=\columnwidth]{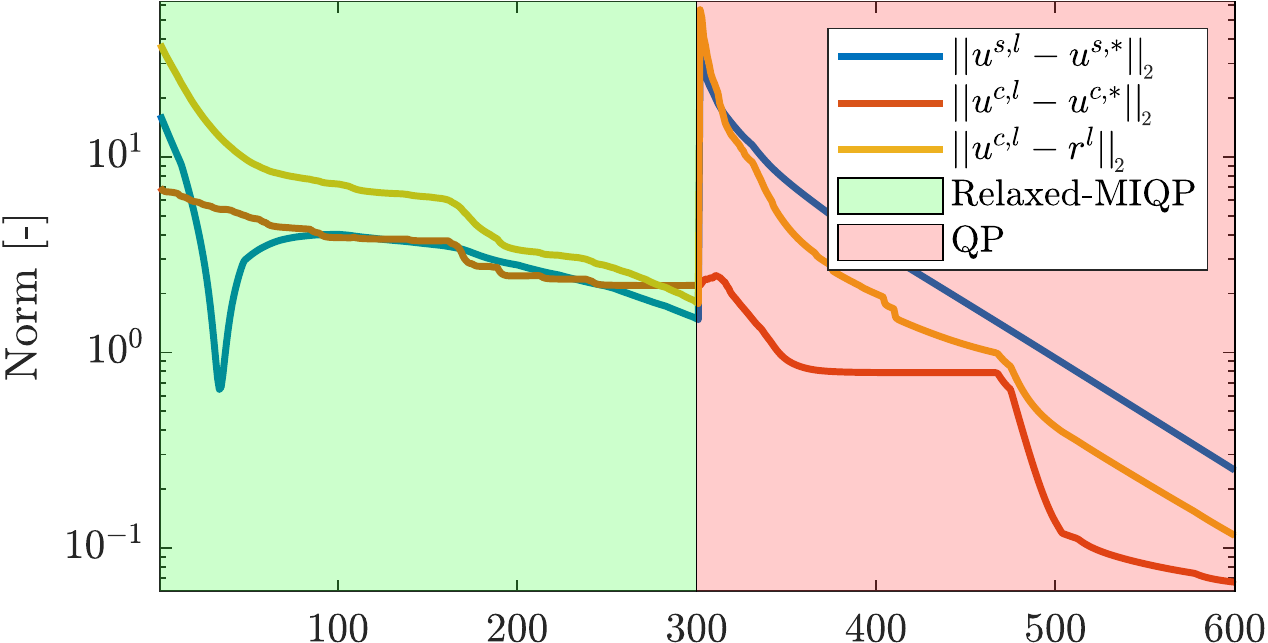}
         \caption{~}
     \end{subfigure}
     \begin{subfigure}[b]{\columnwidth}
         \centering
         \includegraphics[width=\columnwidth]{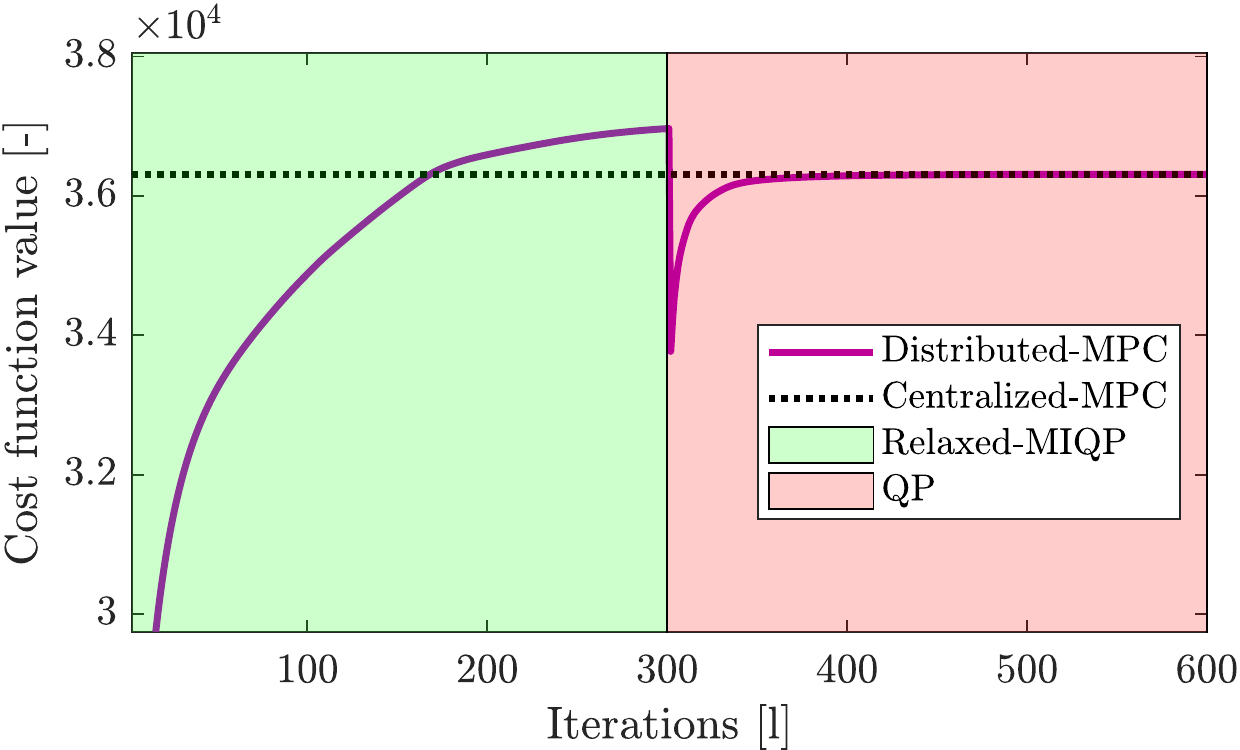}
         \caption{~}
     \end{subfigure}
\caption{Execution of Algorithm~\ref{Algo1} during a time-step of the experiment presented in Section~\ref{experiment}. (a) The difference between the solution of the centralized MPC and the solution of the distributed MPC. (b) The cost function of the centralized MPC is in dotted black and the cost function of the distributed MPC is in solid purple. In both figures, the two stages of the algorithm are differentiated: the semi-relaxed MIQP in the green region and the QP in the red region. The stopping criteria was achieved in 3.38s on a personal computer, an MSI GP62MVR 7RFX Leopard Pro with a 2.5 Ghz Intel i5 7th Core CPU with 8GB of RAM.}
\label{Convergence_cost}
\end{figure}

Algorithm~\ref{Algo1} summarizes the distributed control policy. 
Figure~\ref{Convergence_cost} (a) shows the evolution of the difference between the solution of the centralized MPC and the solution of the distributed MPC at a time-step of the experiment presented in Section~\ref{experiment}. 
Figure~\ref{Convergence_cost} (b) shows the evolution of the cost function of the distributed MPC compared to the cost function of the centralized MPC when  Algorithm~\ref{Algo1} is employed during the same time-step. 
Figure~\ref{Convergence_cost} also displays the two stages of the algorithm where a semi-relaxed MIQP is solved initially to find binary variables heuristically, and subsequently a QP is solved to obtain the optimal solution. 
During step 1, the cost function of the semi-relaxed MIQP converges to a sub-optimal cost with tolerance $\epsilon_{tol,r} = 5 \cdot 10^{-3}$ at iteration 300. 
Given the sub-optimal solution, the semi-relaxed MIQP becomes a QP by fixing the binary variables to 1 if their relaxed counterpart exceeds the threshold $z_{\text{bound}} = 0.5$, and to 0 otherwise. 
If this estimation procedure produces the correct value of the binary variables, the cost function and the solution of the QP respectively converges to the optimal cost and optimal solution with $\epsilon_{tol} = 5 \cdot 10^{-4}$, which is demonstrated in Figure~\ref{Convergence_cost}. 
Note that at the beginning of each stage, the agents start from a feasible solution for their local constraints which is not necessarily feasible for the global optimization problem. 
Then, as they proceed iteratively, they pay for this infeasibility (via the dual variable), and ultimately reach the optimal global (and therefore feasible) solution.
\begin{algorithm}[t]
\DontPrintSemicolon
\SetAlgoLined
\Init{{\normalfont Define the constants $\epsilon_{\text{tol,r}}$, $l_{\text{max,r}}$, $N_{\text{relax}}$, $N$, $\kappa$, $\epsilon_{\text{tol}}$, $l_{\text{max}}$ and $z_{\text{bound}}$.}}{}
\KwIn{At iteration $k$, get current states and disturbance forecasts. Set $p = 0$, $l = 0$ and $c_{\text{tmp}} = \infty$.}
\KwOut{Apply $u^\alpha_{i}(0)$ to the agents at iteration $k$.}
\While{$\epsilon_\text{r} < \epsilon_{\text{tol,r}} \; \text{and} \; l < l_{\text{max,r}} $}{

Apply relaxation to~\eqref{P3} by setting $z_i^{\alpha}(k) \in [0, 1]$, for iteration $k = N_{\text{relax}}, \dots, N$.\;
Solve minimisation problems in~\eqref{P3} to determine the terms $r^{}_i$ and $u^{c_{\text{net}}}_{j}$. \; Update:  $\mathrm{p}_i=\mathrm{p}_i+\kappa(\mathrm{u}_{i}^{c_{\text{net}}}-\mathrm{r}_i)$. \;
Compute: $c = \sum\limits_{i=1}^{N_\text{S}} f^s_i(v^s_i) + \sum\limits_{i=1}^{N_\text{C}} f^c_i(v^c_i)$. \;
Set: $\epsilon_r = |c_{\text{tmp}} - c|$,  $c_{\text{tmp}} = c$ and $l$ = $l$ + 1.}
Impose in~\eqref{P3}: $z_i^{\alpha}(k) = \left\{
    \begin{array}{ll}
        0 & \mbox{if } z_i^{\alpha}(k) < z_{\text{bound}} \\
        1 & \mbox{if } z_i^{\alpha}(k) \geq z_{\text{bound}}
    \end{array}
\right.$, $\forall k$. \\
\While{$\epsilon < \epsilon_{\text{tol}} \; \text{and} \; l < l_{\text{max}} $}{
Solve minimisation problems in~\eqref{P3} to determine the terms $r^{}_i$ and $u^{c_{\text{net}}}_{j}$. \; Update:  $\mathrm{p}_i=\mathrm{p}_i+\kappa(\mathrm{u}_{i}^{c_{\text{net}}}-\mathrm{r}_i)$. \;
Compute: $c = \sum\limits_{i=1}^{N_\text{S}} f^s_i(v^s_i) + \sum\limits_{i=1}^{N_\text{C}} f^c_i(v^c_i)$. \;
Set: $\epsilon = |c_{\text{tmp}} - c|$,  $c_{\text{tmp}} = c$ and $l$ = $l$ + 1.
}

 \caption{Distributed MPC}
 \label{Algo1}
\end{algorithm}

\section{Numerical Study}
\label{Results}

Numerical simulations of the proposed control approach detailed in Section~\ref{methodology} are presented here. These simulations are performed using a model derived from historical data of the building presented in Section~\ref{case_study}. A simulation of a multi-agent environment comparable to a city district is performed in Section~\ref{multi_agents}. A large-scale environment simulation to analyse the computational complexity of the control structures is discussed in Section~\ref{large-scale}.

\subsection{Case Study: NEST}
\label{case_study}

\begin{figure*}[ht!]
	\centering
	\begin{minipage}{.352\textwidth}
		\centering
		\subfloat[]{\includegraphics[width=\textwidth,keepaspectratio]{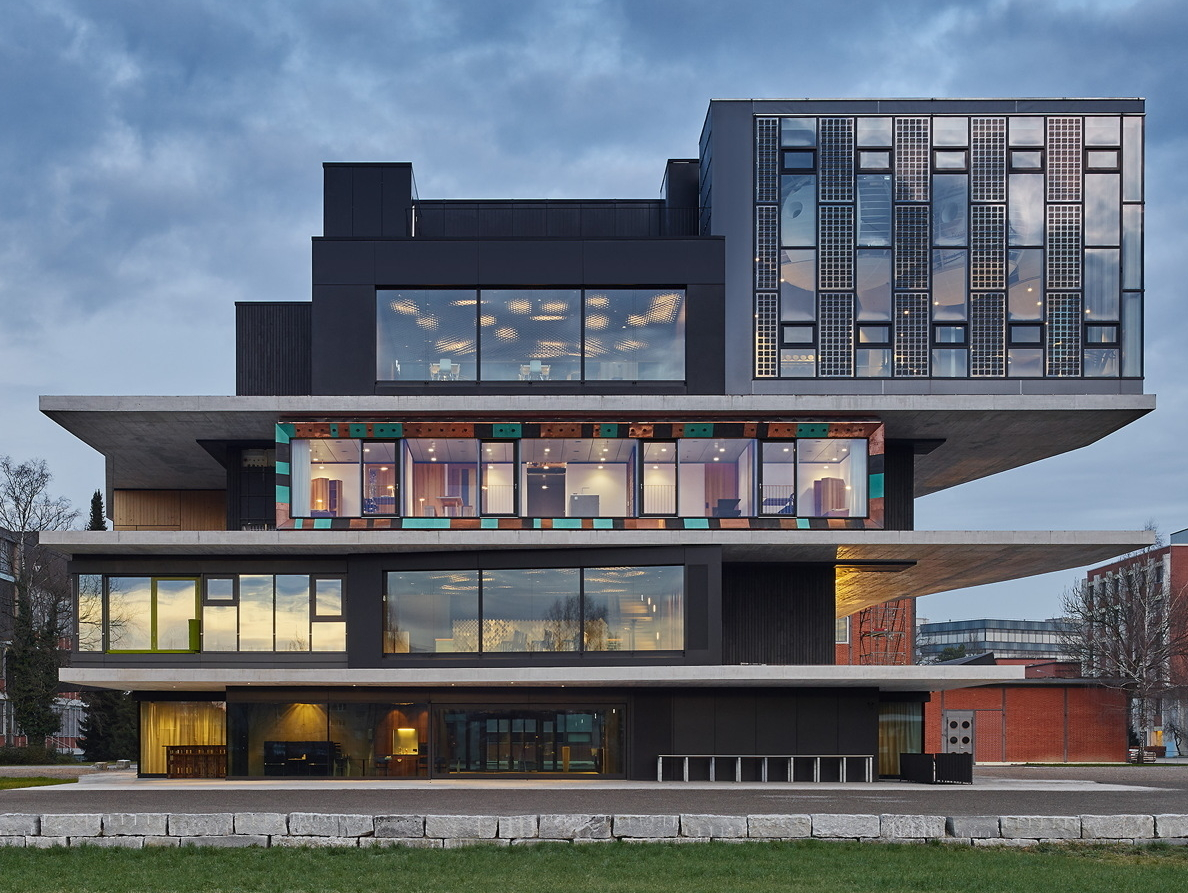}}
	\end{minipage}%
	\hspace{0.16cm}
	\begin{minipage}{0.30\textwidth}
		\centering
		\subfloat[]{\includegraphics[width=\textwidth,keepaspectratio]{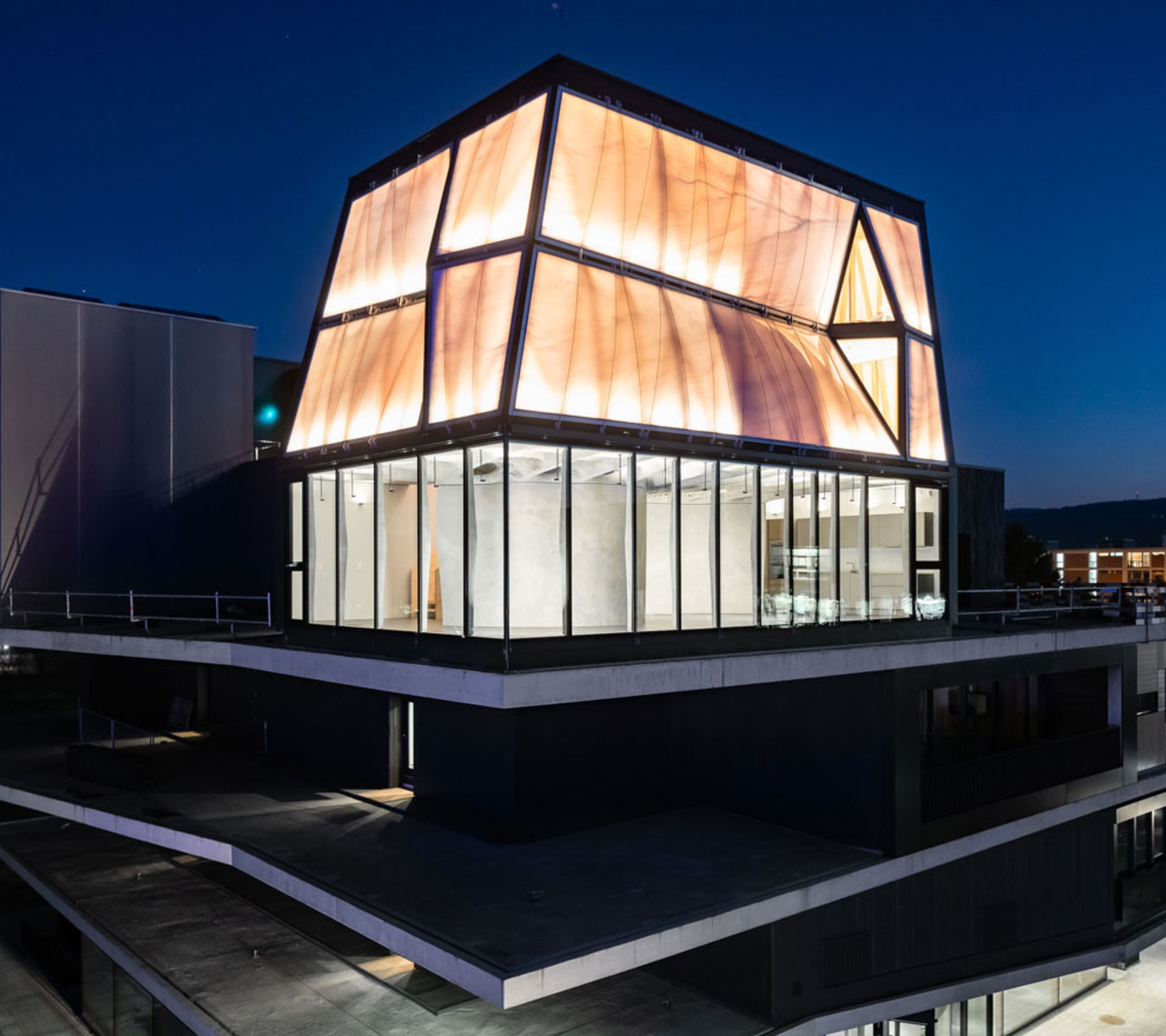}}
	\end{minipage}%
	\hspace{0.07cm}
	\begin{minipage}{0.1145\textwidth}
		\centering
		\subfloat[]{\includegraphics[width=\textwidth,keepaspectratio]{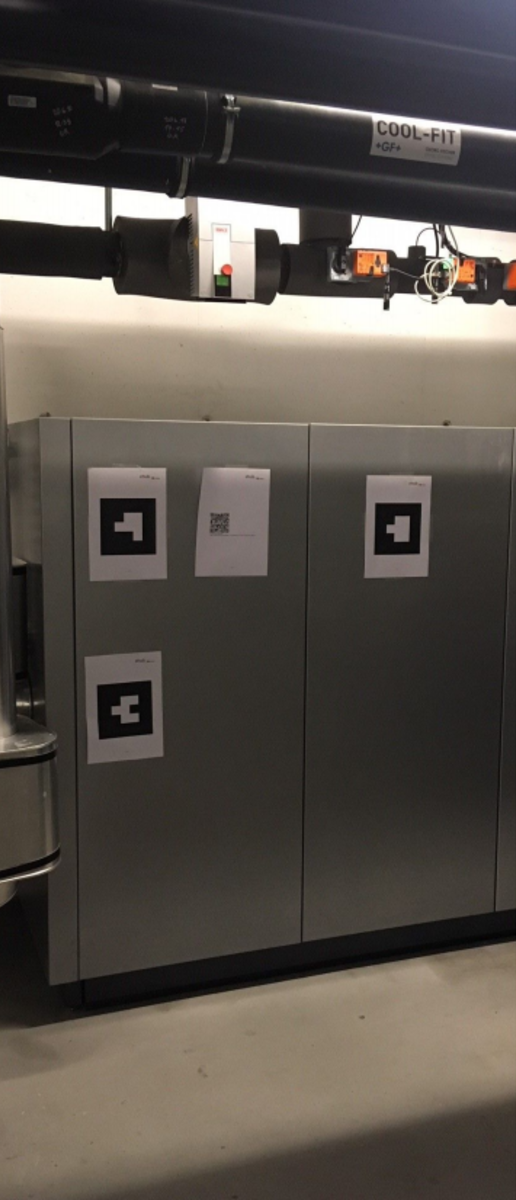}}
	\end{minipage}%
	\hspace{0.07cm}
	\begin{minipage}{0.177\textwidth}
	\centering
	\hspace{0.095cm}
	\subfloat[]{\includegraphics[width=\textwidth,keepaspectratio]{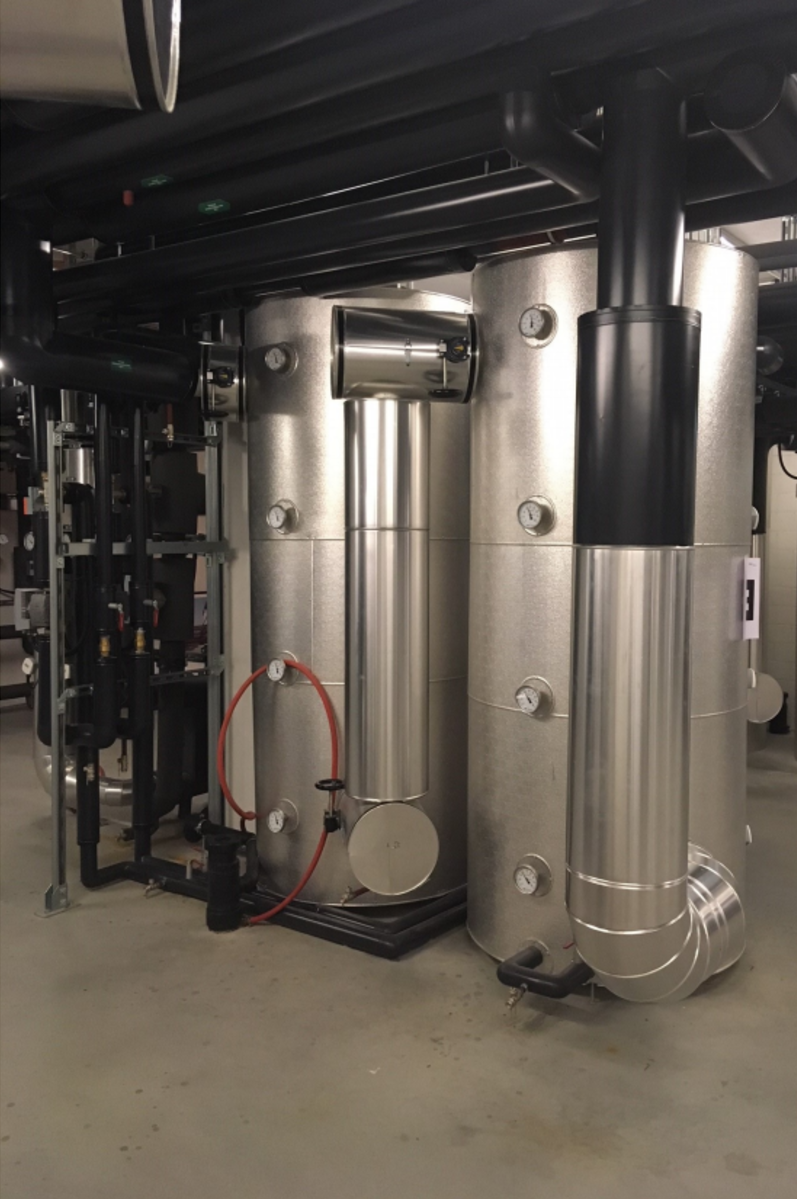}}
	\end{minipage}
	\caption{Different units making up the NEST building (a) Exterior view of NEST. The UMAR unit is the apartment located on the second floor. Copyright: Wojzech Zawarski, Zooey Braun (b) Exterior view of DFAB at night. Copyright: Roman Keller, Tom Mundy and Andrei Jipa (c) the Heat pump and (d) the water storage of NEST.  Copyright: Felix Bünning.}
	\label{NEST}
\end{figure*}

\begin{figure}[b!]

\centering
\begin{subfigure}[b]{\columnwidth}
  \centering
  \includegraphics[width=\columnwidth]{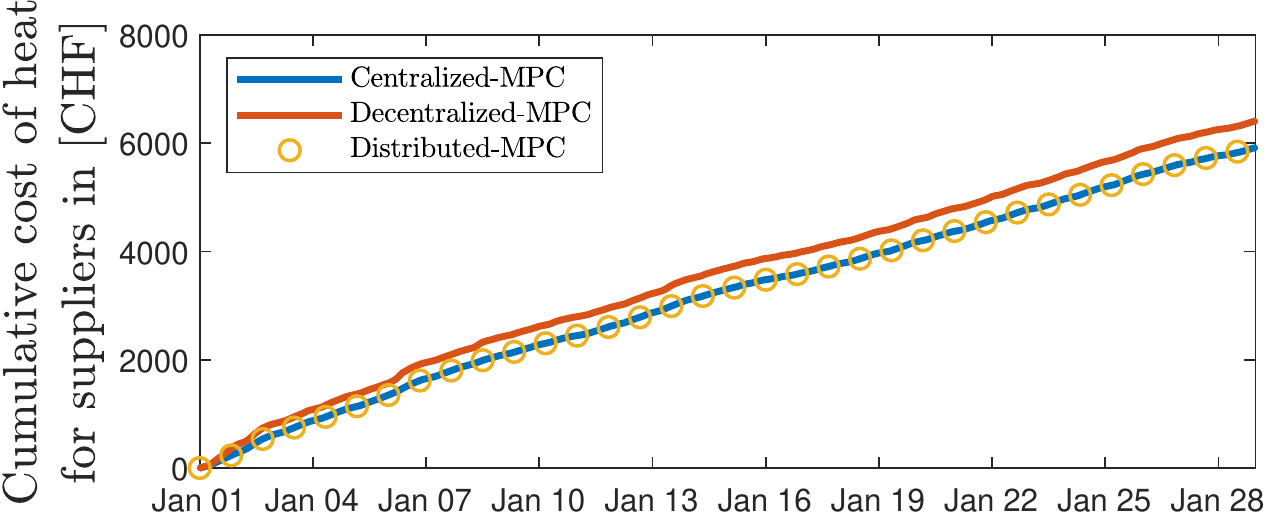}
  \caption{~}
\end{subfigure}
\begin{subfigure}[b]{\columnwidth}
  \centering
  \includegraphics[width=\columnwidth]{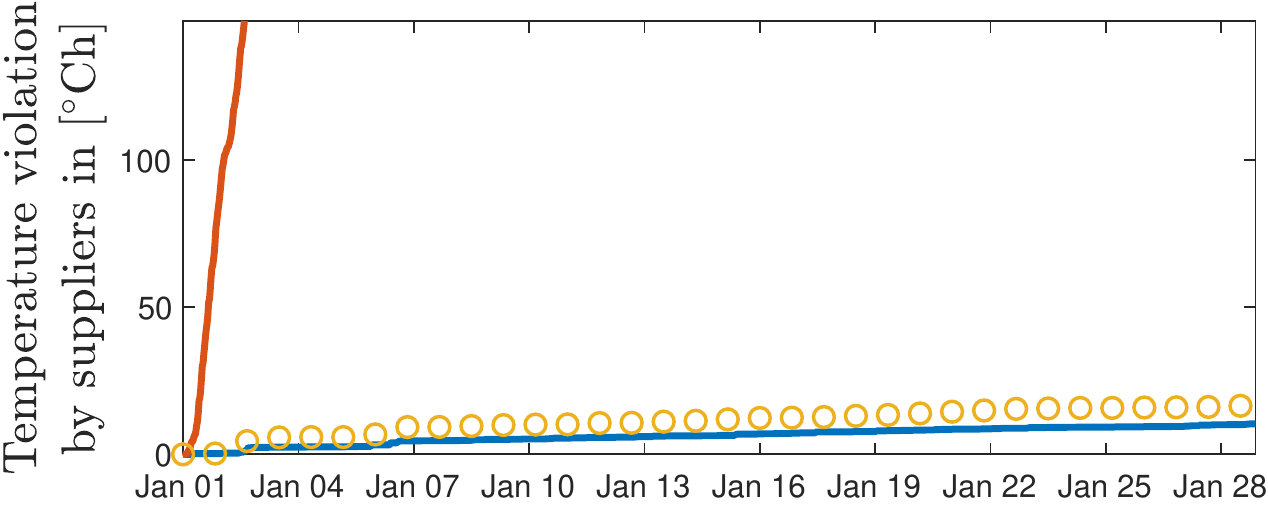}
  \caption{~}
\end{subfigure}
\begin{subfigure}[b]{\columnwidth}
  \centering
  \includegraphics[width=\columnwidth]{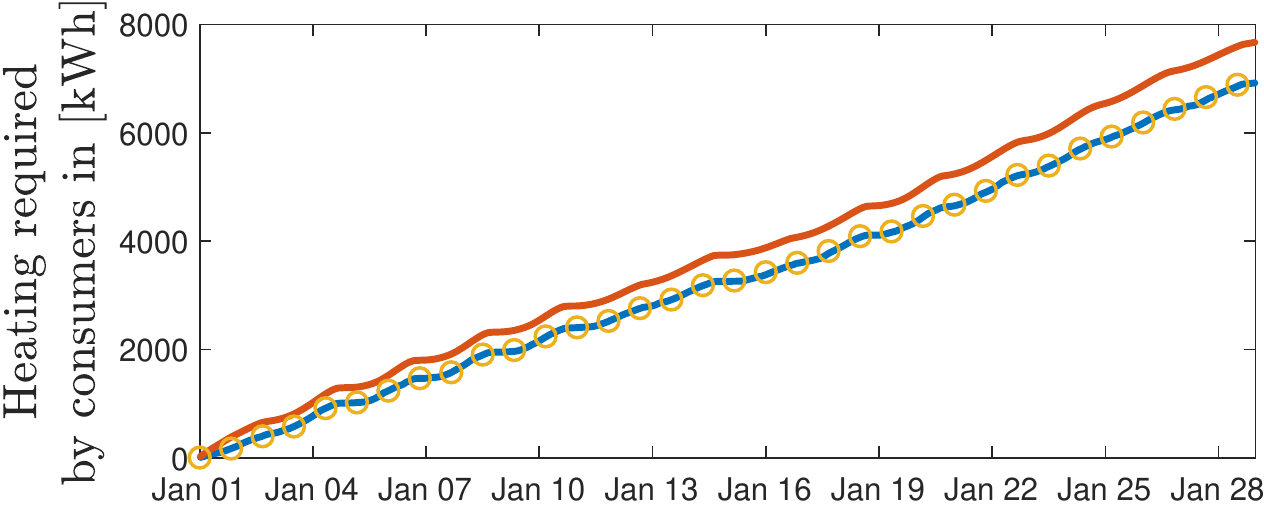}
  \caption{~}
\end{subfigure}
\begin{subfigure}[b]{\columnwidth}
  \centering
  \includegraphics[width=\columnwidth]{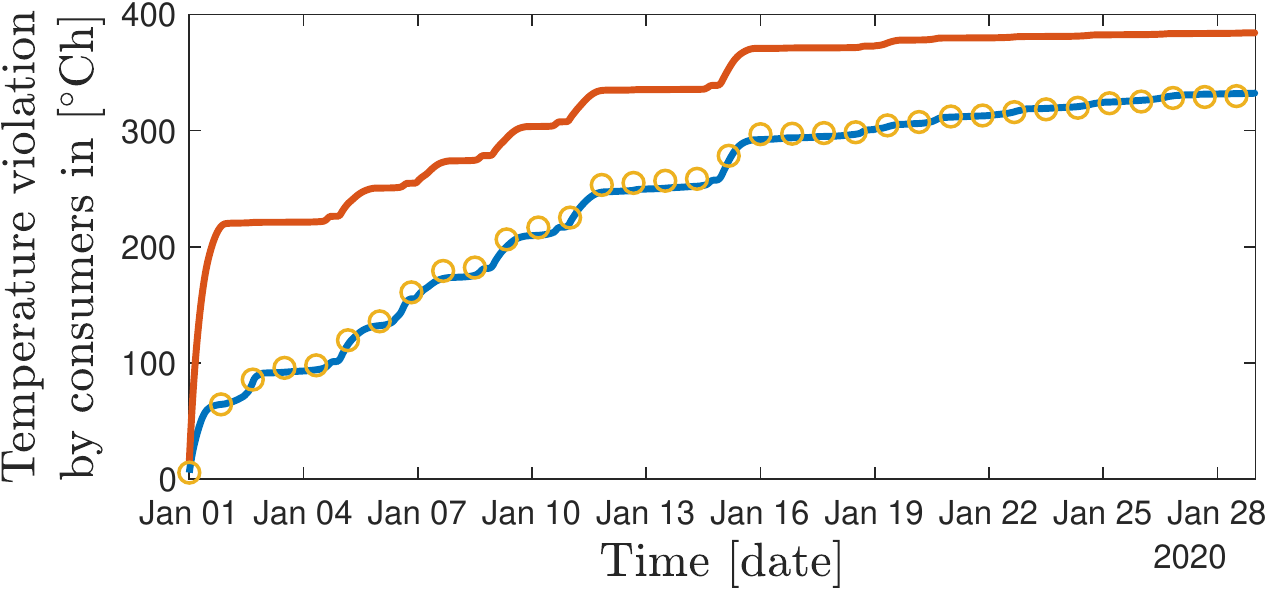}
  \caption{~}
\end{subfigure}

\caption{Performance of centralized, decentralized and distributed control schemes. (a) the cumulative heat entering in the apartments. (b) the cumulative temperature constraint violation in the tanks. (c) the cumulative heat entering the tanks. (d) the cumulative temperature constraint violation in the rooms.}
\label{Consumer_supplier_res}
\end{figure}

\begin{figure*}[p]

\centering

\begin{subfigure}[b]{2\columnwidth}
  \centering
  \includegraphics[width=0.8\textwidth]{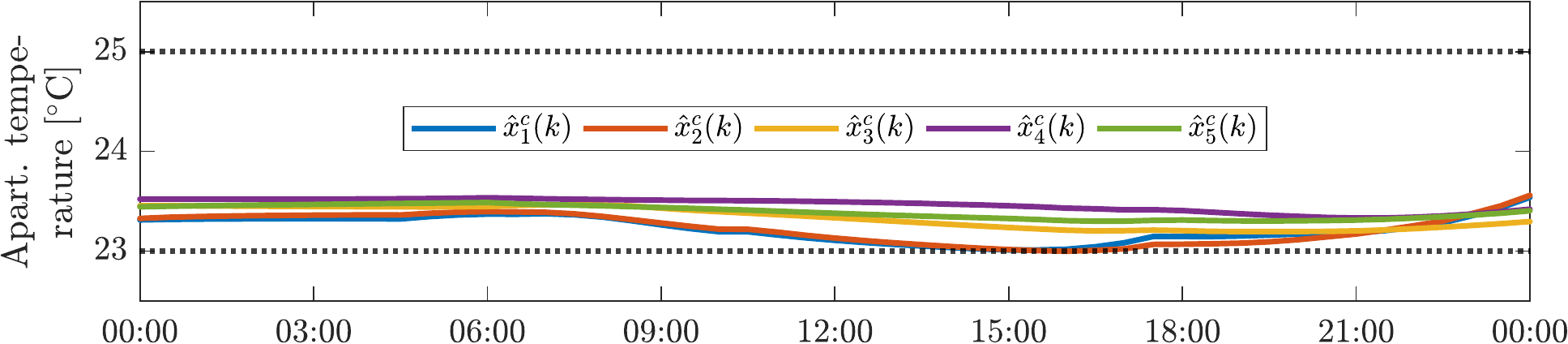}
  \caption{~}
\end{subfigure}

\begin{subfigure}[b]{2\columnwidth}
  \centering
  \includegraphics[width=0.8\textwidth]{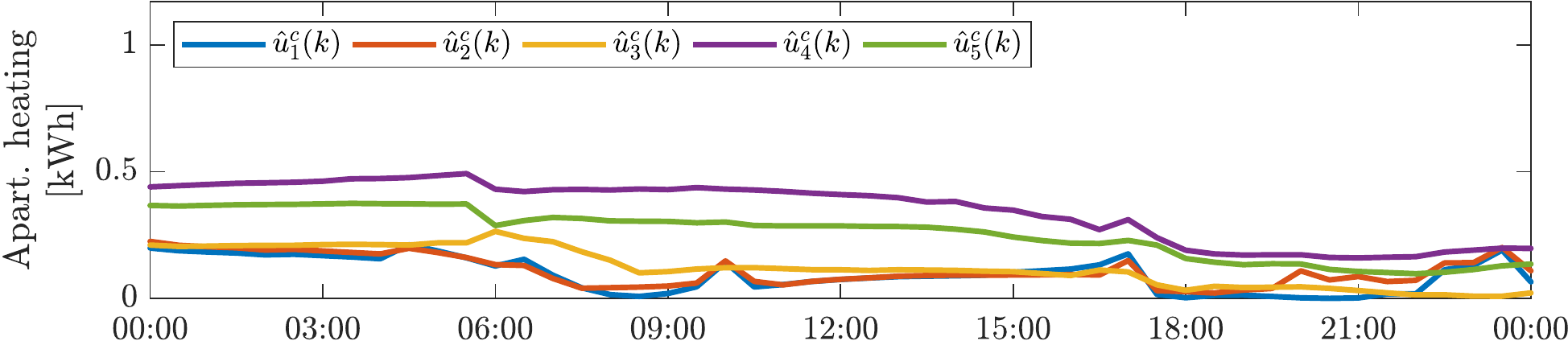}
  \caption{~}
\end{subfigure}

\begin{subfigure}[b]{2\columnwidth}
  \centering
  \includegraphics[width=0.8\textwidth]{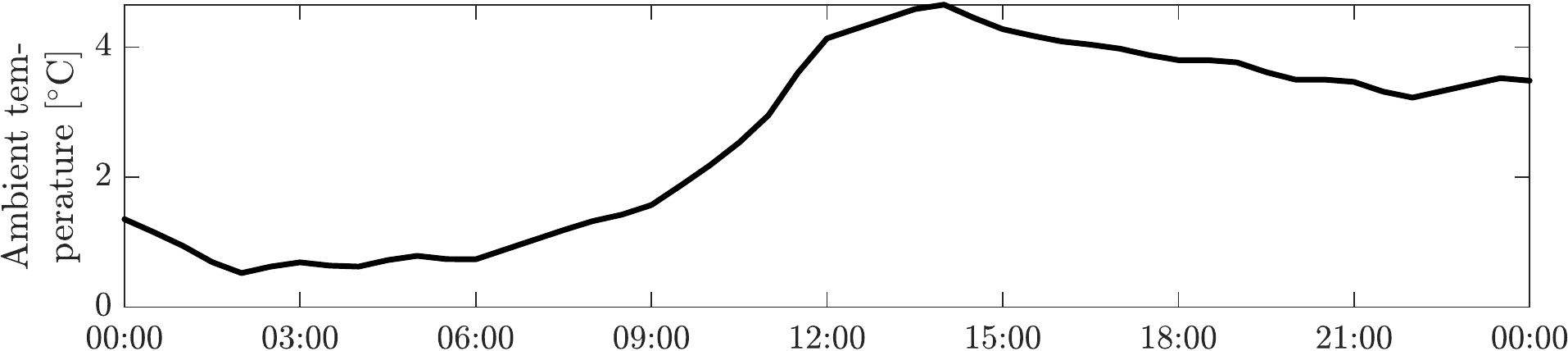}
  \caption{~}
\end{subfigure}

\begin{subfigure}[b]{2\columnwidth}
  \centering
  \includegraphics[width=0.8\textwidth]{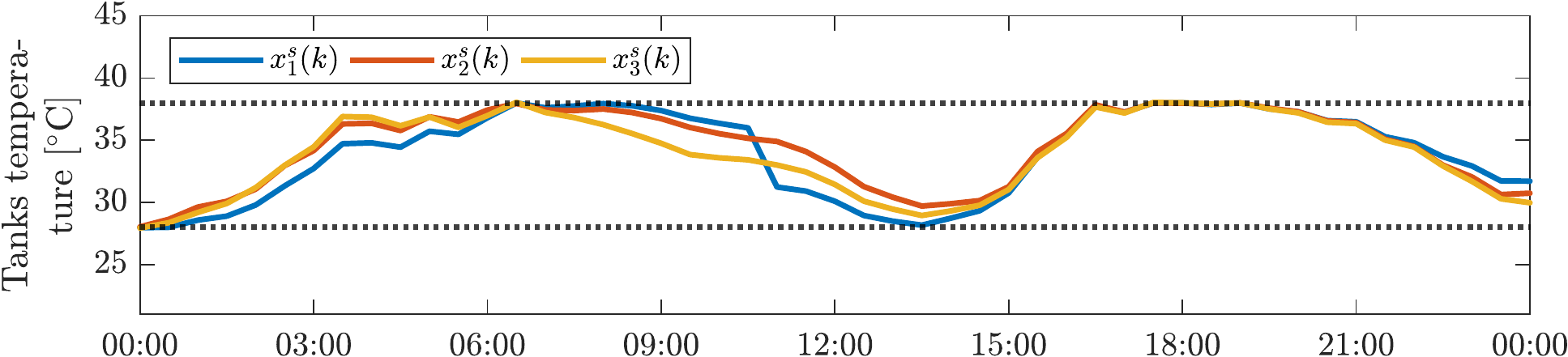}
  \caption{~}
\end{subfigure}

\begin{subfigure}[b]{2\columnwidth}
  \centering
  \includegraphics[width=0.8\textwidth]{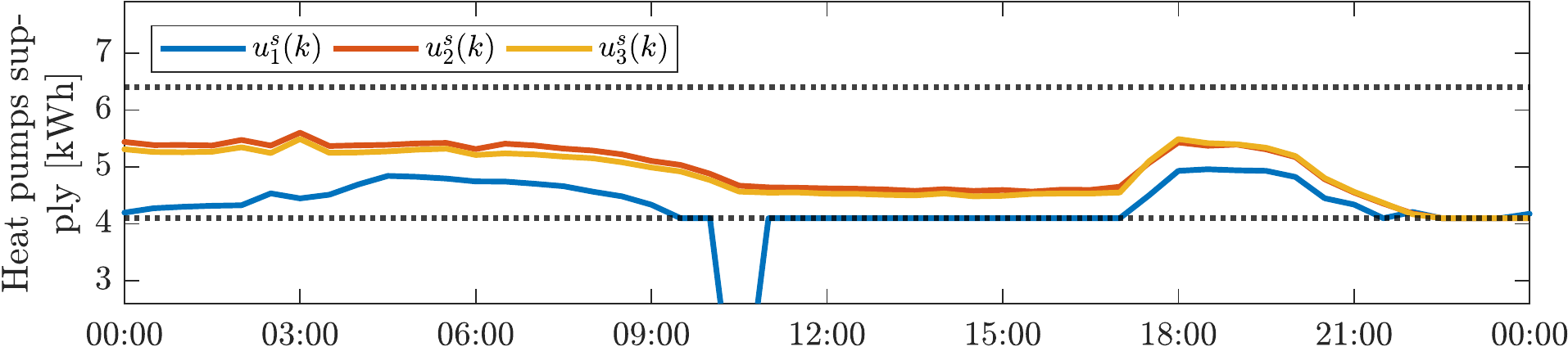}
  \caption{~}
\end{subfigure}

\begin{subfigure}[b]{2\columnwidth}
  \centering
  \includegraphics[width=0.8\textwidth]{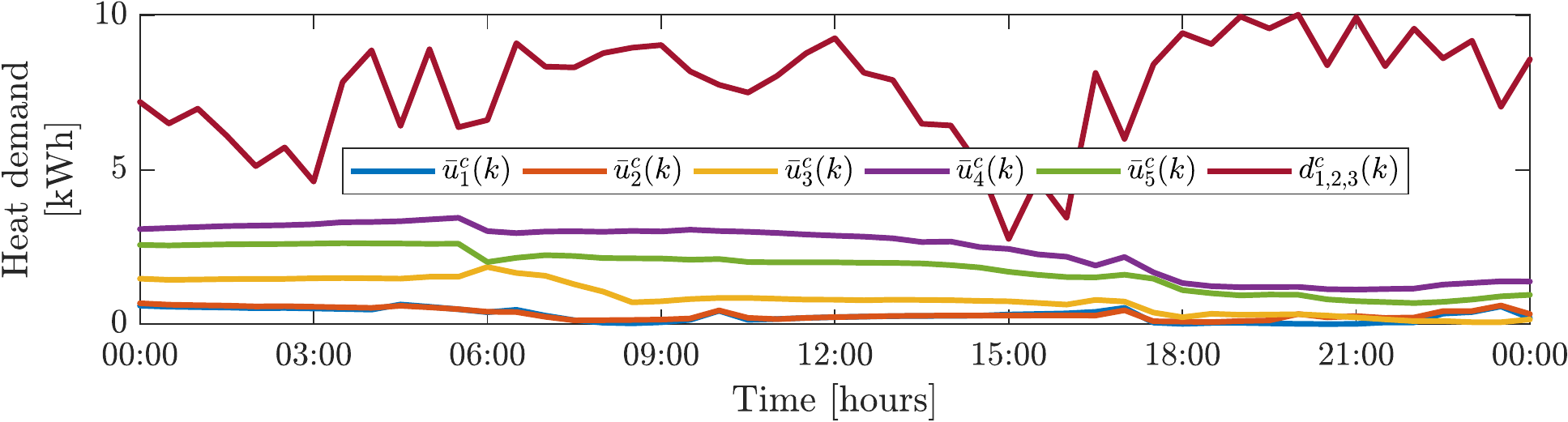}
  \caption{~}
\end{subfigure}

\caption{Numerical experiment results. (a) Average temperature of the rooms of each apartment, black dotted lines indicate temperature constraints (b) Average room heat supply in each apartment (c) Ambient temperature outside the NEST building (d) Average tank temperatures, black dotted lines indicate temperature constraints (e) Heat supply of each heat pump, black dotted lines indicate input constraints (f) Heating demand of the uncontrolled buildings and the apartments.}
\label{num_experiment}
\end{figure*}

The NEST (Next Evolution of Sustainable Building Technology) building~\cite{Richner_Heer_Largo_Marchesi_Zimmermann_2017} is an energy hub demonstrator at Empa in D\"{u}bendorf, Switzerland. The aim of the demonstrator is to test new technologies, materials, and systems in terms of their impact on energy management in buildings. A picture of the facility is shown in Figure~\ref{NEST} (a). The building hosts a wide variety of technologies that convert and store energy. It also comprises various units with different use cases (residential, offices, meeting rooms) that can be temporarily installed in the NEST core structure. All units have individual heating and cooling system and an individual control system. Thus, in our paradigm, NEST can be viewed as a simulator of the interaction between buildings and energy hubs.

In the context of this study, we employ three agents from the demonstrator: the Urban Mining and Recycling (UMAR) unit and the Digital Fabrication unit (DFAB) as two consumers, and the medium temperature grid with a water buffer storage supplied by a heat pump as one supplier.

The UMAR and DFAB units are both apartments in the NEST building. UMAR is shown in Figure~\ref{NEST} (a) (also in caption) and DFAB is shown in Figure~\ref{NEST} (b). The original purpose of UMAR is to demonstrate the uses of fully reusable, recyclable, or compostable resources in construction~\cite{brunner2011urban}. The DFAB unit is distinctive in that it was not only digitally designed and planned but also built using predominantly digital processes, both on-site and off-site~\cite{menna2020opportunities}. Both units comprise seven rooms each. The units are equipped with heating systems that take their energy from the medium-temperature grid (with a supply temperature between 28 \si{\celsius} and 38 \si{\celsius}) of the NEST building via heat exchangers. The heat is then sent to the rooms through pipes and distributed by ceiling heating panels in UMAR and a floor heating system in DFAB. The heat transferred to each room can be estimated by combining the supply valve position for each room and the total energy consumption of the unit. Note that in this study, control is only available for three rooms in UMAR: 272, 273, and 274, which are the bedrooms and the living room of the unit. The heat supply in rooms 272 and 274 is constrained to 0.6 \si{\kilo\watt} while room 273 is constrained to 1.8 \si{\kilo\watt}. For DFAB, the heat supply in rooms 371, 472, 474, 476 571, 573, and 574 are respectively constrained to 0.83, 0.55, 1.36, 1.39, 1.14, 1.35 and  0.75 \si{\kilo\watt}. Each room is equipped with a temperature sensor. Forecasts for the disturbances, which are ambient temperature and solar irradiation, are available from MeteoSwiss~\cite{SwissFedMeteo}. The mathematical representation of the units was obtained by gray-box modelling parameters estimated from historical data captured from sensors in the units. 

The supplier agent of the medium temperature heating system comprises a ground-source heat pump and a water storage tank. The devices are shown in Figure~\ref{NEST} (c) and (d). The heat pump draws cold water from the bottom of the storage, warms it up by transferring heat from the refrigerant to the water inside the condenser, and feeds it back into the top of the storage tank. 
The heat demand of the units is met with individual pumps drawing warm water from the top of the storage tank and passing it through heat exchangers, where the heat is transferred to the units' heating systems.  
The average conversion efficiency between electrical energy and high-temperature thermal energy in the heat pump is described by the coefficient of performance $\alpha_{COP} = 3.53$. 
The electrical capacity of the heat pump is between 8.2~\si{\kilo\watt} and 12.8~\si{\kilo\watt}. The mathematical representation of the storage and the heat pump is based on first-principles models established with simple thermal heat transfer equations, which have been validated via  experiment~\cite{Bnning2020MachineLA}. 
As storage tanks are industrial products, we assume that the parameters can be obtained from the manufacturer or easily be determined with high accuracy.

\subsection{Multi-Agent Simulation}
\label{multi_agents}
In this section, we present the results of numerical experiments in order to evaluate the performance of the different control schemes presented in Sections~\ref{centralizedMPC}-\ref{DistributedMPC}. 
The objective is to simulate an environment comparable to a small city or a district, i.e., control of multiple hubs and dwellings over a long duration. Using historical data and identified models of units of the NEST building, the simulation was conducted using disturbance data from January 1-28, 2021.

\setlength{\tabcolsep}{5pt} 
\renewcommand{\arraystretch}{1.12}
\begin{table*}[!b]
    \caption{Comparison table between centralized, decentralized and distributed control approach.}
    \label{Comparison_table}
\centering
\begin{tabular}{c|c|c|c|c|}
\cline{2-5}
                                        & \begin{tabular}[c]{@{}c@{}}Room heating \\ in {[}kWh{]}\end{tabular} & \begin{tabular}[c]{@{}c@{}}
                                        Room comfort zone \\ violation in {[}$\degree$Ch{]}\end{tabular} & \begin{tabular}[c]{@{}c@{}}Tank heat \\ supply in {[}kWh{]}\end{tabular} & \begin{tabular}[c]{@{}c@{}}Tank temperature constraint \\ violation in {[}$\degree$Ch{]}\end{tabular} \\ \hline
\multicolumn{1}{|c|}{Centralized-MPC}   & 6920.6                                                               & 332.2                                                                                                          & 21285 (5917.1 CHF)           & 10.4                                                                                                           \\ \hline
\multicolumn{1}{|c|}{Decentralized-MPC} & 7672.2                                                               & 384.1                                                                                                          & 21824 (6408.5 CHF)           & 746.8                                                                                                         \\ \hline
\multicolumn{1}{|c|}{Distributed-MPC}   & 6949.5                                                               & 330.4                                                                                                          & 21314 (5937.6 CHF)            & 16.7                                                                                                           \\ \hline
\end{tabular}
\end{table*}

The simulation environment comprises eight agents. There are three suppliers (three tanks with heat pumps and boilers, representative of that used in NEST), i.e. $N_{\text{S}} = N_{\text{HP}} = N_{\text{HB}} = 3$. Each tank is connected with an uncontrolled building, i.e. $N_{\text{NC}} = 3$.  There are five consumers (two UMAR-like apartments and three DFAB-like apartments), i.e. $N_{\text{C}} = 8$. The first supplier is connected to the first and second consumers, the second supplier is connected to the second, third, and fourth consumers, and the third supplier is connected to the fourth and fifth consumers. The characteristics of the agents (states, inputs, constraints) were kept the same as described in Section~\ref{case_study}, while dynamic matrices of the agents were perturbed to make them  non-homogeneous. For the three suppliers, the coefficient of performance $\alpha_{COP}$ occurring in the input matrices $B^s_i$ were modified from the original value and respectively set to 3.7095, 3.6728, and 3.5367. Additionally, as a backup, an electric boiler with a coefficient of performance of 1 and an electrical capacity of 0 kW to 20 kW was connected to each tank. For the five consumers, the dynamic matrices were obtained for different indexes $i$ by taking the original matrix components equal to their nominal values plus a random term normally distributed with a mean of zero and a standard deviation of $\sigma_i$. The term $\sigma_i$ was defined as the standard deviation of the set containing the components of matrix $A^s_i$. The initial condition of the simulation, i.e., the starting temperatures of the tanks and the rooms of the apartments, were randomly chosen outside the corresponding operational constraints within a margin of 2\si{\celsius}. The external disturbances (heating demand of uncontrolled buildings, ambient temperature, and solar irradiance) were obtained from actual measurements from the NEST building during the specified simulation period.

The following simulation parameters were used in our evaluation of centralized MPC, decentralized MPC, and distributed MPC. 
The time-step between two adjacent control steps was set to 30 minutes. 
The prediction horizon was set to $N = 24$, which corresponds to 12 hours. The upper and lower comfort constraints of all rooms are 23 \si{\celsius} and 25 \si{\celsius} respectively. 
To simulate varying electricity prices, the coefficient of the suppliers' weighting matrix $Q^s_i$ in the cost function were defined based on the local scheduled electricity tariff, i.e.~17.07 cents/\si{\kilo\watt\hour} for off-peak between 10 p.m.~and 6 a.m.,~and 28.06 cents/\si{\kilo\watt\hour} for on-peak during the rest of the day. 
The input weighting matrices of the consumers $Q^c_i$ were set to identity.  The weighting matrices $R^s_i$ and $R^c_i$ were set to $R$ = $200 \cdot I_{N} \otimes I$. 

In the case of distributed MPC, the following additional parameters were used. 
The convergence tolerances were set to $\epsilon_{\text{tol,r}} = 5 \cdot 10^{-3}$ and $\epsilon_{\text{tol}} = 5 \cdot 10^{-4}$. 
The minimum number of iterations was set to 150, and the maximum number of iterations for the relaxed-MIQP was set to $l_{\text{max,r}} = 300$, while for the QP it was set to $l_{\text{max}} = 850$. 
This limitation is only used here to simulate possible operating delays. 
The relaxation horizon $N_{\text{relax}}$ was set to 12, the decision boundary $z_{\text{bound}}$ was set to 0.5, and the dual variable step size was set to $\kappa = 0.15$. 
Finally, the weighing factor $\rho_{ij}$ was tuned experimentally and set to $0.08$ for each pair of connected suppliers $i$ and consumers $j$.

Table~\ref{Comparison_table} presents a comparison of the results of each of the tested MPC controllers. 
It shows the cumulative heat consumption and constraint violation at the end of the simulation for both the suppliers and the consumers. 
As the centralized controller has complete information about all agents, the solution is considered as the true optimum and serves as a reference for the other controllers. 
It can be seen that the decentralized strategy shows poor results in terms of constraint violations. 
Compared to the centralized strategy, the suppliers have violated the tank temperature constraints by a factor of over 70. 
Moreover, the overall system consumes more energy as the consumers requested a heat surplus of 10.1\% compared to the centralized controller. 
The distributed strategy delivers improved  performance. 
The consumers only need a heat surplus of 0.42\% compared to the central solution while violating the constraints 0.54\% less. 
The results show that the distributed MPC scheme has a comparable performance to the centralized MPC, while both methods significantly outperform the decentralized MPC. 

Figure~\ref{Consumer_supplier_res} shows the cumulative costs for the heating suppliers, the integrated temperature constraint violation in the tanks, the cumulative heat entering the apartments, and the integrated temperature constraint violation in the rooms. 
We can observe that the centralized and distributed MPC have very similar performance  in terms of heat consumption and constraint violations. The decentralized MPC shifts gradually from the optimal performance over time. 
This is due to the fact that the coupling between the agents is not considered by the individual optimization problems. For instance, without knowing the consumers' heat demand, the suppliers do not predictively adapt their heat supply from the tanks. 
As the capacity of the heat pumps is limited, this results in a violation of the tanks' state constraints, as can be seen in Figure~\ref{Consumer_supplier_res} (b). 
Similarly, without knowledge of the suppliers' maximum capacity, the consumers consume too much heat compared to the minimum required, as pictured in Figure~\ref{Consumer_supplier_res} (c). 
Note that the high violations of consumers comfort constraints visible in figure~\ref{Consumer_supplier_res} (d) is a result of initial conditions outside of the operational constraints at the beginning and an undersized heating system present in the real UMAR and DFAB units.

Figure~\ref{num_experiment} shows detailed trajectories of the distributed controller for a single day of the numerical experiment.
The historical data used for the disturbances is from January 25, 2021. Figure~\ref{num_experiment} (a) shows the average temperature of the rooms of each apartment in colored bold lines along with the temperature constraints in dotted black lines. It can be seen that all temperatures remained in the comfort zone during the day. The average temperatures are shown for better visibility; meanwhile, we note that the individual temperatures also stayed within constraints. It can be seen that most of the temperatures stay closer to the lower bound during the whole experiment. During times of high ambient temperatures (see (c)), between 12:00 and 15:00, some of the temperatures are close to the constraint, as the heat pump can be expected to have enough capacity reserves during these times. The optimization-based strategy determines that it is an unnecessary use of energy to have temperatures too high above the lower limit.

Figure~\ref{num_experiment} (b) shows the average heat input to the rooms of each apartment, while Figure~\ref{num_experiment} (c) depicts the measured ambient temperature outside of NEST. We can see that the energy supply to the rooms coincides well with the evolution of the outside temperature. Indeed, in every apartment, the heat input is reduced approximately 50\% from beginning to end of the day when the outside temperature doubled in the same period.

Figure~\ref{num_experiment} (d) shows the average storage temperature of each tank in solid colored lines and the temperature constraints in dotted black lines. Figure~\ref{num_experiment} (e) shows the heat generated by each heat pump in solid colored lines and the corresponding operational constraints in dotted black.  Figure~\ref{num_experiment} (f) depicts the (uncontrolled) actual heating demand of the NEST building in solid red and the (controlled) heating demand of the apartments in solid colored lines. In Figure~\ref{num_experiment} (d), we can see that the temperatures stayed within the constraints from the beginning to the end of the experiment. In Figure~\ref{num_experiment} (e), we can see that the inputs of the suppliers are, to a large extent, determined by the NEST heat demand. Nevertheless, the apartments' heat consumption also affects the input of the suppliers, as the heat pumps inject more heat in the tanks at the beginning of the day than at the end.

\subsection{Numerical Analysis for Large-Scale Scenarios}
\label{large-scale}

According to Section~\ref{multi_agents}, the centralized and distributed controllers have  comparable  performance in terms of energy consumption and constraint violations. In contrast, the decentralized controller has poor performance, and thus, it is not suitable to be implemented in practice. However, the decentralized approach is less computationally demanding than the other two, because the optimization calculation is partitioned between the agents and executed simultaneously. The distributed controller has the same calculation configuration, but requires iterating between the agents to reach  a consensus on the shared optimal variables. For implementation on real systems, where the optimization time is constrained, it is necessary to investigate the computational demand of the centralized and distributed controllers. In order to address this question, large-scale simulations, i.e.,~involving a large number of agents, are performed in the sequel.

Starting from an environment with two agents (one supplier and one consumer) up to an environment with $n$ agents, we calculated the computational time necessary for each controller to find its solution for a single time-step. In each environment, the suppliers correspond to the one described in Section~\ref{case_study}, i.e.~a water storage supplied by a heat pump and an electric boiler supplying an uncontrolled building heating demand. The consumers are UMAR-like apartments and were generated using the same perturbation method as described in~\ref{multi_agents}. The electrical capacity of the suppliers' heat pumps was scaled according to the number of connected consumers. The rest of the characteristics specific to the simulation are the same as those defined in Section~\ref{multi_agents}. The controllers are configured precisely as in Section~\ref{multi_agents}.

For a single  time-step to find optimal control inputs, the computational time taken by each controller type was calculated as follows. For the centralized controller, the computation time per time-step was considered as the time required for one MPC to solve for the control action of all agents. For the distributed controller, we compiled a list containing all individual computational times required by each local MPC to solve for the control action of its agent. The computational time per time-step was then defined to be the maximum time in this list. A summation of the computational time of each iteration was performed until convergence is reached to obtain the total time of one time-step of control actions.

In order to simulate a realistic environment, we defined scenarios in which the diversity of agents followed specific  rules. The quantity of suppliers with respect to consumers was defined by a ratio and the ratio was always less than 1, i.e., there are more consumers than suppliers.  In addition, consumers were evenly distributed among suppliers, where one consumer is always shared between two suppliers to have a coupling between the networks. In this fashion, we carried out a series of scenarios in which the ratio between the number of suppliers and consumers varied for a fixed number of agents. For each scenario, starting from two agents, consumers and suppliers were added to the environment up to a maximum of $n = 200$ agents.

\begin{figure}[H]
\centering
\includegraphics[width=\columnwidth]{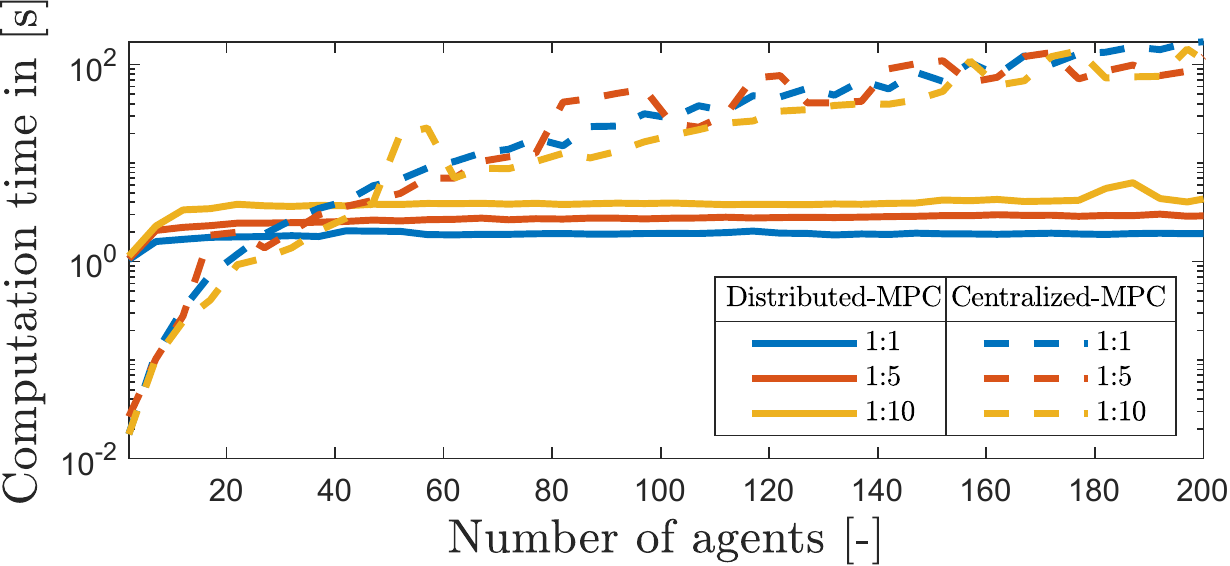}
\caption{Computational analysis. Total computational time necessary to calculate the control action for the centralized and distributed control schemes with respect to the number of agents for different ratios between the number of suppliers and consumers. The simulation was performed using the Euler cluster at ETH~\cite{EulerETH:2021}.}
\label{time_constant}
\end{figure}

Figure~\ref{time_constant} shows the computational demand of the centralized and distributed controllers for different scenarios. For the distributed controller, the computational time always ends up becoming constant for systems with many agents, while it consistently increases for the centralized controller. The intersection point for the number of agents where the distributed controller becomes less computational demanding than the centralized controller increases with the ratio, i.e., 22 agents for 1:1, 32 agents for 1:5, and 37 agents for 1:10. Indeed, in the case where the number of consumers managed by a single supplier increases significantly, the computational demand of the supplier agent increases as the substitute variable $\bar{r}$ increases in dimension. In the case of the centralized controller, this variable $\bar{r}$ increases in size as the number of agents, whereas in the case of the distributed controller, for a fixed ratio of consumers for a single supplier, the local vectors $r_i$ do not increase in size, which explains the stabilization of the computation time.

The results indicate that implementing a centralized control law for large-scale environments is difficult in practice because of high computational requirements. Another disadvantage of a central solution is its lack of resilience: damage to the central controller will cause the failure of the entire energy management system. This is not the case for a distributed control system. Note that the communication time between agents in decentralized and distributed control is not modeled here, although we can reasonably expect it to be small compared to the optimization solving time.

Note also that the centralized controller uses the Gurobi solver while the solver for the distributed controller uses a mixture of Gurobi and a custom solver: the MPC of each individual agent solves its local optimization problem using Gurobi, but the external unit solves the dual problem using a simple heuristic sub-gradient method. The comparison of the absolute calculation time between the two controllers is therefore irrelevant. In the case where the dual problem is also solved with a potentially faster commercial solver, the computation time of the distributed controller should be even lower than what is presented in this study.

\section{Experimental Validation}
\label{experiment}

\begin{figure*}[!t]

\centering

\begin{subfigure}[b]{2\columnwidth}
  \centering
  \includegraphics[width=0.8\textwidth]{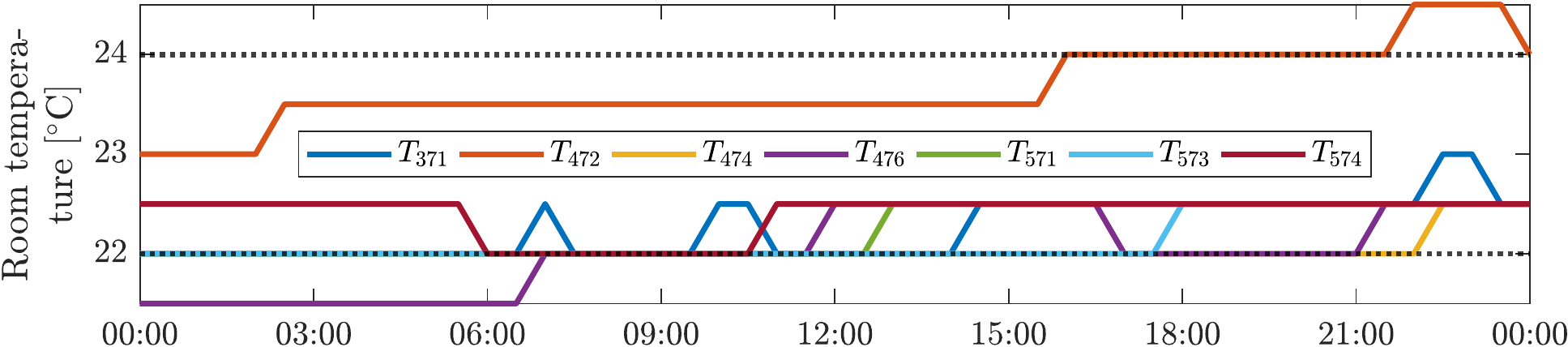}
  \caption{~}
\end{subfigure}

\begin{subfigure}[b]{2\columnwidth}
  \centering
  \includegraphics[width=0.8\textwidth]{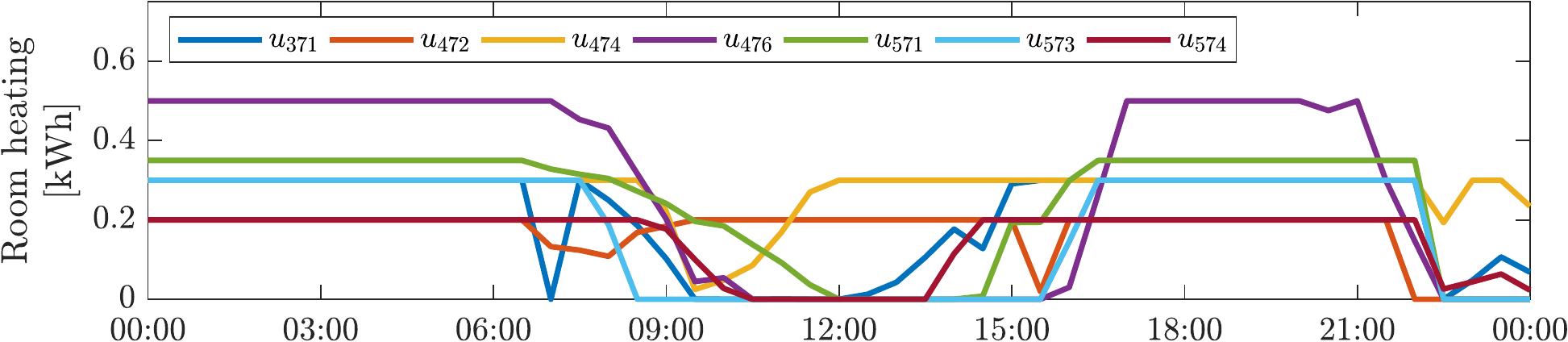}
  \caption{~}
\end{subfigure}

\begin{subfigure}[b]{2\columnwidth}
  \centering
  \includegraphics[width=0.8\textwidth]{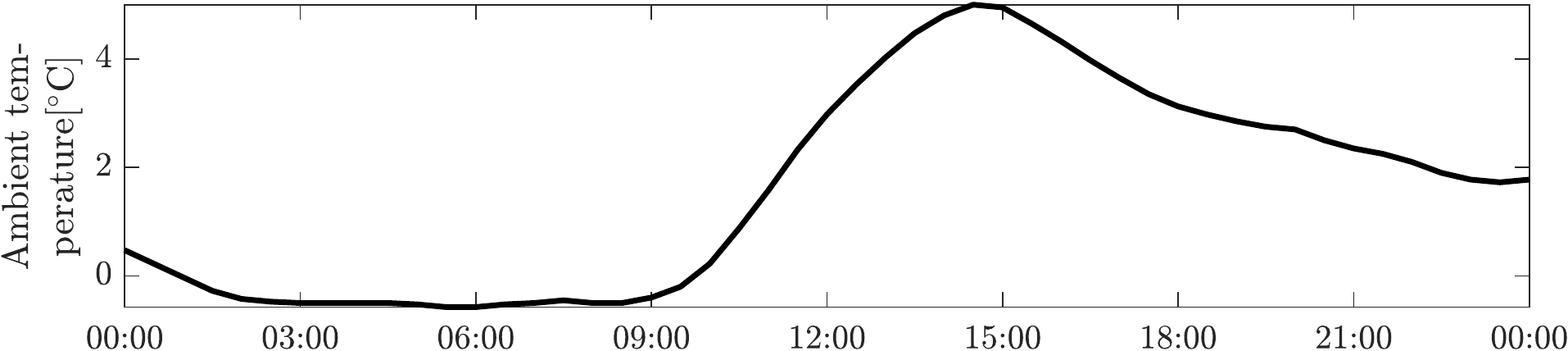}
  \caption{~}
\end{subfigure}

\begin{subfigure}[b]{2\columnwidth}
  \centering
  \includegraphics[width=0.8\textwidth]{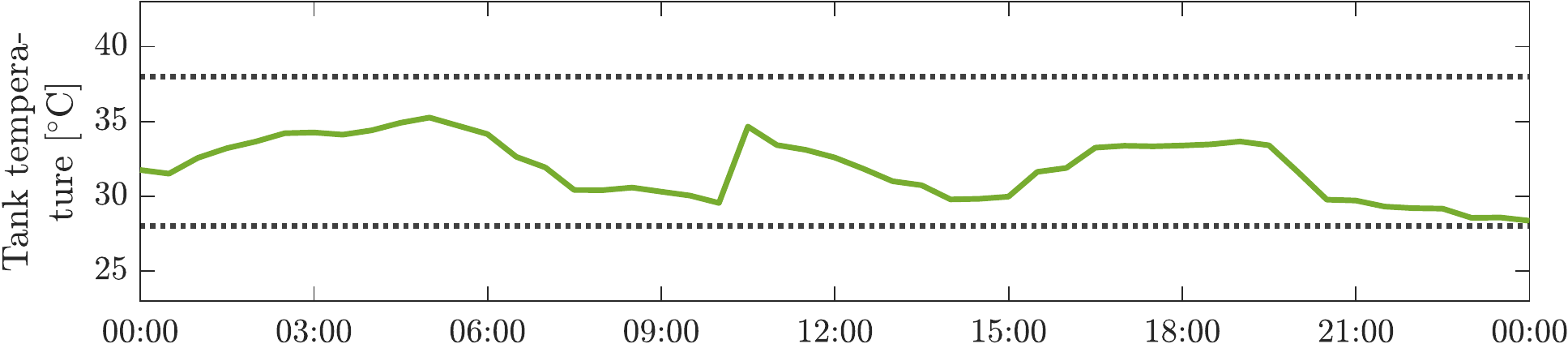}
  \caption{~}
\end{subfigure}

\begin{subfigure}[b]{2\columnwidth}
  \centering
  \includegraphics[width=0.8\textwidth]{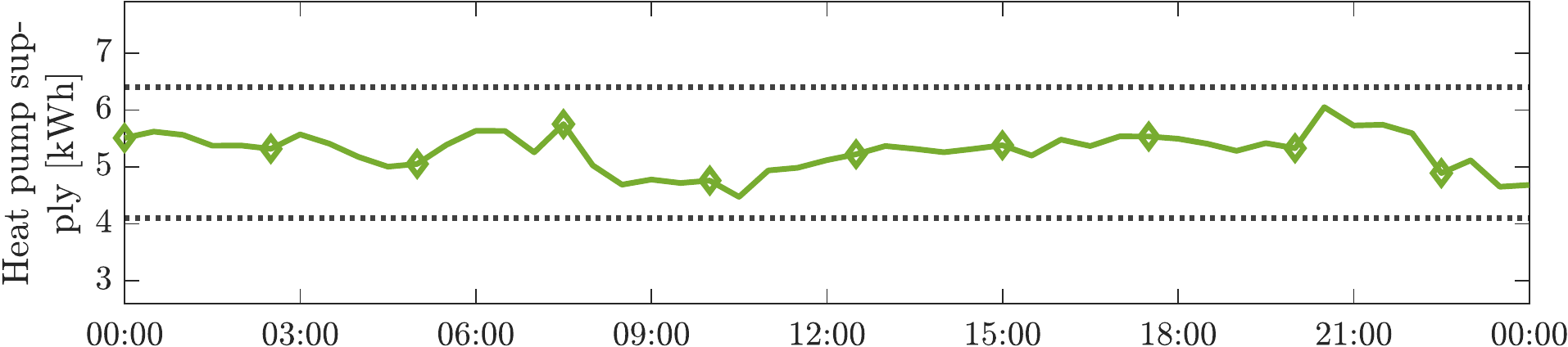}
  \caption{~}
\end{subfigure}

\begin{subfigure}[b]{2\columnwidth}
  \centering
  \includegraphics[width=0.8\textwidth]{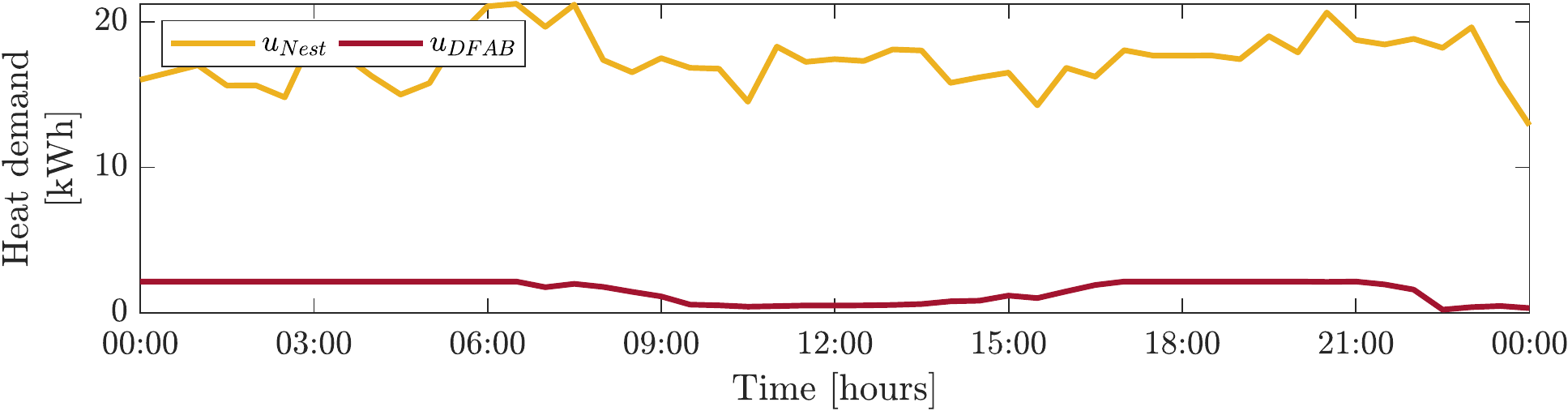}
  \caption{~}
\end{subfigure}

\caption{Experimental results. (a) Temperature in each room of DFAB unit, black dotted lines feature temperature constraints (b) Heat supply in each room of DFBA unit (c) Ambient temperature outside NEST building (d) Average tank temperature, black dotted lines feature temperature constraints (e) Heat input delivered by the heat pump, black dotted lines feature input constraints (f) Heating demand of NEST building and DFAB unit.}
\label{experiemtn}
\end{figure*}

To evaluate the performance of the proposed approach in a representative application, the distributed controller was tested during an experiment involving two agents in the NEST building: one supplier (the NEST medium-temperature heat pump and a water buffer storage) and one consumer (the DFAB unit). 
The experiment was conducted over the period of 24 hours starting at 00:00 on December 15, 2020.

Most of the controller parameters remain as in Section~\ref{multi_agents}, however some modifications were made to the specific configuration due to the presence of occupants. The comfort constraints were set to 22 \si{\celsius} and 24 \si{\celsius} for the rooms of the DFAB apartment. Moreover, the coefficients of the supplier's input weighting matrix $Q^s_1$ were set to identity. Finally, no electric boiler was used during the experiment, as the backup system installed in NEST is operated by a standard controller, in case that an experiment causes the heat pump to fail.

The controller and the related optimization schemes were implemented in \texttt{MATLAB} and solved with \texttt{Gurobi}. In each time-step, the calculation of the MPC control inputs was started three minutes before they were applied to the agents. This time period is sufficiently long to carry out the optimization, i.e.~for the agents to converge to a solution. The system states were estimated by extrapolating measured states of the previous time-step and the measured states at the beginning of the optimization. The communication between the agents was implemented in \texttt{Python 3}. A \texttt{Python} OPC-UA client was used for the communication with sensors and actuators of the agents. Additional \texttt{Python} 3 scripts were used for the forecast of the heating demand of NEST. These are out of the scope of this study, but are discussed in detail in~\cite{BUNNING2020109821}. 

The results of the one-day experiment are shown in Figure~\ref{experiemtn}. Figure~\ref{experiemtn} (a) shows the temperature of the rooms in DFAB in colored bold lines, along with the temperature constraints in dotted black lines. It can be seen that all temperatures remained within the comfort constraints, except for the temperature in room 476 at the start of the experiment and the temperature in room 472 at the end of the experiment. For room 476, this is simply due to the given initial condition. After 6 a.m.~the temperature reaches the lower comfort  constraint and stays above until the end. In the case of room 472, the temperature violates the upper comfort constraint for two hours at the end of the experiment. This is likely an effect of the coarse granularity of the temperature sensor (0.5 \si{\celsius}) and could be mitigated with the help of a state estimator. At the end of the experiment, the temperature satisfies the constraint.

Figure~\ref{experiemtn} (b) shows the heat input to each room of DFAB, and Figure~\ref{experiemtn} (c) shows the measured ambient temperature outside of NEST during the experiment. The energy supplied to the rooms coincides well with the evolution of the ambient temperature. At the beginning of the day, the exterior temperature is low, and consequently, all inputs are at their maximum. In the middle of the day, the ambient temperature rises, and most of the rooms significantly reduce their heat consumption. Moreover, it can be seen that the controller expects the ambient temperature to rise and thus stops heating early: while the ambient temperature only starts to rise significantly at 09:00, the controller already reduces heating in most rooms between 07:00 and 08:00. At the end of the day, the ambient temperature drops again, resulting in an increase of heating.

Figure~\ref{experiemtn} (d) shows the average storage temperature in solid green and the temperature constraints in dotted black. Figure~\ref{experiemtn} (e) shows the heat generated by the heat pump in solid green with diamond markers and the operational range of the pump in dotted black. Finally, Figure~\ref{experiemtn} (f) depicts the real heating demand of the NEST building in solid yellow and the heating consumption of DFAB in solid red. In Figure~\ref{experiemtn} (d), we can see that the temperature stays well within the constraints throughout the entire experiment. Note that the temperature of the tank is only marginally affected by the energy demand of DFAB but rather driven by the demand of the rest of the NEST building. Indeed, plot (f) shows that the energy demand of NEST is much larger than the demand of the DFAB unit.

Unfortunately, a direct comparison between the distributed, centralized and decentralized MPC approaches cannot be made in real experiments, as the experimental conditions are not repeatable. However, the case study indicates that practical implementation of the distributed controller performs satisfactorily in an occupied budding application.
\section{Conclusion}
\label{sec:conclusion}

In this study, we have developed an MPC control structure for the management of energy in an environment where both energy hubs and buildings are considered as controlled entities. Three different approaches have been studied: centralized, decentralized, and distributed. Extensive numerical experiments modeling a building-scale energy hub system showed that the distributed approach was the most appropriate solution. In the considered environment, the method managed to offer good performance with low computational load by exploiting the output coupling between agents through virtual shared prices. Furthermore, the method avoids sharing local constraints and states, which reduces the need for agents to share potentially private information. An experimental implementation was performed on a building and energy hub located in D\"ubendorf, Z\"urich, to demonstrate the practical feasibility and the effectiveness of the method. The results obtained were satisfactory as it was capable of maintaining the room comfort constraints by taking into account external disturbances and optimizing energy consumption.

Future work will focus on testing the experimental implementation on longer periods and varying configurations to demonstrate the controller's robustness. Extension of the simulation environment to other technologies, including cooling systems, photovoltaics, and batteries, will be investigated. Another direction of research will be the investigation of alternative optimization methods to improve the performance and the speed of the distributed controller.
\section*{Acknowledgement}

This work was supported by the  Swiss  Competence  Centers for  Energy  Research  FEEB\&D project and the  ETH  Foundation.  The authors would like to thank the Urban Energy Systems Laboratory Empa, particularly Sascha Stoller and Reto Fricker, for helping with the coordination and access to the NEST demonstrator.


\end{document}